\documentclass[12pt]{article}
\usepackage{amsmath}
\usepackage{amssymb}
\usepackage{graphicx}

\usepackage{cite}
\usepackage[british]{babel}

\newtheorem{theo}{Theorem}[section]
\newtheorem{lm}{Lemma}[section]
\newtheorem{df}{Definition}[section]

\newtheorem{cor}{Corollary}[section]

\newtheorem{proposition}{Proposition}[section]

\allowdisplaybreaks

\numberwithin{equation}{section}

\def\R{{\mathbb R}}
\def\Z{{\mathbb Z}}
\def\N{{\mathbb N}}
\def\S{{\mathcal S}}
\def\Exp{{\mathbf E}}
\def\E{{\mathbf E}}
\def\Pr{{\mathbf P}}
\def\1{{\mathbf 1}}

\def\eps{\varepsilon}

\newcommand{\tg}{{\mathop{\rm tg}}}
\renewcommand{\arctan}{{\mathop{\rm arctg}}}

\newcommand{\proof}{\noindent {\bf Proof.}~}

\def\qed{\hfill\rule{.2cm}{.2cm}}

\newcommand{\ud}{{\mathrm d}}
\newcommand{\F}{{\mathcal F}}
\newcommand{\DD}{{\mathcal D}}
\newcommand{\fr}{\partial\DD}

\title{Asymptotic behaviour of randomly reflecting billiards
in unbounded tubular domains}

\author{M.V.~Menshikov\thanks{Partially supported by FAPESP (07/50459--9)}$^{~,1}$ \and 
M.~Vachkovskaia\thanks{Partially supported by CNPq (304561/2006--1, 471925/2006--3) and FAPESP (thematic grant 04/07276--02)}$^{~,2}$ \and  A.R.~Wade\thanks{Partially supported by
         the Heilbronn Institute for Mathematical Research}$^{~,3}$}

\begin{document}

\maketitle
{\footnotesize
\noindent
$^1$ Department of Mathematical Sciences,
University of Durham,
South Road, Durham DH1 3LE, UK.\\
E-mail: \texttt{Mikhail.Menshikov@durham.ac.uk}

\noindent
$^2$ Department of Statistics,
Institute of Mathematics, Statistics and Scientific Computation,
University of Campinas--UNICAMP,
P.O. Box 6065, CEP 13083--970, Campinas, SP, Brazil.\\
E-mail: \texttt{marinav@ime.unicamp.br}

\noindent
$^3$ Department of Mathematics,
University of Bristol,
 University Walk, Bristol BS8 1TW, UK. \\
 Tel: +44 (0) 117 954 6954;
 Fax: +44 (0) 117 928 7999;
 E-mail: \texttt{Andrew.Wade@bris.ac.uk}}

\begin{abstract}
We study stochastic billiards in
 infinite planar domains   with curvilinear boundaries: that is,
 piecewise deterministic motion
 with randomness introduced via
 random reflections at the domain boundary.
 Physical motivation for the process
 originates with ideal gas models in the Knudsen regime,
 with particles reflecting 
 off microscopically rough
 surfaces.
  We classify
the process into  recurrent and transient cases. 
We also
give almost-sure results on the long-term behaviour of the location
of the particle, including a super-diffusive rate of escape in the transient case. 
A key step in obtaining our results
is to relate our process to an instance 
of a one-dimensional stochastic process
 with
asymptotically zero drift, for which we
prove some new almost-sure bounds of independent
interest.  We obtain some of these bounds
via an application of general semimartingale
criteria, also 
of
some independent interest.
\end{abstract}



\section{Introduction}
We consider stochastic billiards (Knudsen random walk) in two-dimensional infinite domains (`tubes')
of the form 
$\DD = \{ (x,y) \in \R^2 : x>A, |y| < g(x) \}$ where $A\ge 1$ and  $g:[1,\infty) \to (0,\infty)$ is a  monotone smooth function.
Generally speaking, billiards are dynamical systems describing the motion of a particle in a region
with reflection rules at the boundary:
such systems have been extensively studied
in the mathematical and physical literature. When the reflection rule is randomized, we have a stochastic billiard
process. 
The stochastic models have received much less attention;
invariant distributions for
stochastic billiards in general (mostly bounded) domains were studied in \cite{CPSV,evans}. In the present
paper we study 
the qualitative (recurrence or transience)
and quantitative (almost-sure bounds)
behaviour of
stochastic billiards in {\em unbounded} domains in the plane.

Physical motivation for the billiards
model comes from the dynamics of
ideal gas models in the so-called
Knudsen regime where intermolecular
interactions are neglected. The random reflection
law is motivated by the fact that
the particle is small and the
surface off which it reflects
has a complicated (rough) microscopic structure. The behaviour of
such Knudsen flows is of interest in several areas of physics, chemistry,
and technology. For physical background,
see for instance \cite{cer,knu}. For more
on the motivation of the model in the present paper,
see \cite{CPSV} and references therein. We mention some related
models at the end of this section.

Informally, the model that we study can be described as follows. A particle moves ballistically
inside the domain $\DD$ at constant 
velocity 
(of unit magnitude, say),
 and each time that
  the particle hits the boundary, it is reflected at a
  random angle $\alpha \in (-\pi/2,\pi/2)$ to
 the
 inwards-pointing normal at
  the boundary
 curve $\fr$, where $\alpha$ does not depend on the initial direction of the particle.
So, inside the domain the motion is deterministic: the stochasticity is introduced
via the distribution $\alpha$ of the random reflections. In this paper,
we take the distribution of $\alpha$ to be symmetric about $0$.

We consider two types of problem: (i) the continuous-time
motion of the particle
 in the tube; and (ii)
the discrete-time embedded process obtained
by observing the instances of
collisions on the boundary (in other words,
unit time elapses between reflections).
The embedded process
is, from our point of view,
of interest in its own right (and
behaves very differently to the continuous-time
process), and
is also a vital ingredient
for the study of the continuous-time process.
The asymptotic
phenomena that we study for each process
 are also of two
main types: (i) recurrence or transience,
and existence of moments for recurrence times;
and (ii) almost-sure bounds for the location
of the particle at time $t$, as $t \to \infty$.

What kind of regions $\DD$ are of interest?
In the case of $g(x) \equiv 1$, our model becomes
 symmetrically
reflecting
motion in a strip. Roughly speaking for the moment, the horizontal
motion of the particle is then described by a random
walk with zero drift, and so the model is null-recurrent.
On the other hand, if $g(x) = x$ our tube is
a wedge and, with our reflection rule, 
at any reflection there is 
a positive chance that the particle will head off to infinity
and never return to a bounded set: it is transient. 
Similarly if $g(x) = \beta x$ for $\beta \in (0,1)$,
it is not
difficult to show that the
particle will follow nearby the deterministic
 trajectory which goes to infinity with linear speed.
This
dichotomy motivates the study of tubes with widths that
grow sub-linearly, to probe precisely 
the transition
 between recurrence and transience.
It is also natural to consider the case where the tube
has decreasing width. 

The primary family
of tubes to bear in mind has $g(x) = x^\gamma$
for some $\gamma < 1$, although
we do consider more general
forms for the function $g$. Then there
are two main cases: $\gamma \in (0,1)$ and $\gamma <0$. See Figure \ref{tubes}.

We obtain criteria for recurrence and transience
of our processes. Loosely speaking, these results imply, for example,
that for $\gamma \geq 1/2$ the particle always has transient dynamics,
but can be recurrent
for $\gamma \in (0, 1/2)$ depending on the reflection
distribution. Also, for very shallow tubes such as $g(x)=\log x$,
the introduction of random reflections ensures
that the process is recurrent, while the 
corresponding deterministic evolution along the
normals is transient. (See  Theorem \ref{gen_g}
 below.)

We also obtain almost-sure bounds
for the horizontal position of the
particle,
in both discrete- and
continuous-time settings. Thus we
have information about the asymptotic
speed of the processes. For example,
when $g(x) = x^\gamma$ for $\gamma \in [1/2,1)$,
Theorem \ref{time} below
implies that for all large times $t$ the continuous-time
process is at position $t^{1/(2-\gamma)}$, ignoring
logarithmic terms. In particular, the motion is {\em super-diffusive}.

In the case $\gamma <0$, recurrence is evident. In this case
we obtain a criterion for the finiteness of the mean 
recurrence-time (roughly speaking, `ergodicity'). We also
obtain polynomial estimates for the position of the particle;
the ergodicity has so-called `heavy tails'.

A step in our proofs
will involve
relating the stochastic billiard
process to an instance of the
so-called Lamperti problem
(after \cite{lamp1,lamp3,lamp2}) of a
one-dimensional
stochastic process with asymptotically
zero mean drift. A crucial
 ingredient to our proofs
for the stochastic billiard process
is provided by some new results
for general Lamperti-type
processes. These results (Theorems \ref{prop1}, \ref{Lamp_low},
\ref{prop6}) are of independent interest, and are in some
cases apparently
new even for the nearest-neighbour random walk (birth-and-death chain).
To obtain our almost-sure bounds for Lamperti-type
processes,  we state and prove some general results
on obtaining almost-sure bounds for stochastic processes via semimartingale-type criteria. 

We conclude this section with some remarks upon related models.  
 As previously mentioned,
stochastic billiard models
have received comparatively little attention.
Classical (that is, non-stochastic)
billiard models have been extensively studied, particularly in mathematical
physics, and the literature is vast; see for instance
\cite{tab}. A common approach in the classical setting 
is via
dynamical systems for which the reflection rule
is elastic. In this setting,
we mention that billiards in certain
unbounded domains resembling
the tubes considered here (at least
for $\gamma\leq 0$) have been
studied; see for instance
\cite{eml,lenci1,lenci2,lenci3}
and references therein. 
An infinite-tube billiard model 
with a stochastic
component (cf the $\gamma=0$ case here)
is analyzed
in \cite{bab}.
In the dynamical systems
setting, studying ergodicity is typically
a primary goal. In the stochastic
setting, the results of the present paper
demonstrate a complete
range of behaviours, from positive-recurrence
(roughly speaking, `ergodicity'),
through null-recurrence, to transience.

In the next section we give a more formal
definition of the model
and state our main results. The statement
of our auxiliary results
on the Lamperti problem
is deferred to Section
\ref{lamperti}.

\section{Description of the model and results}

\subsection{Construction}

The rigorous formulation
of the stochastic billiard model
that we consider is essentially
given in \cite{CPSV}. We now describe the
construction, which is modified
slightly from that in \cite{CPSV}
to fit with our context.  Let $A \in [1,\infty)$, which
will be   large but fixed, to be specified later. 
For a monotone function $g:[1,\infty) \to (0,\infty)$,
let
\[ \DD := \DD (g;A) := \{ (x,y) \in \R^2 : x>A, |y| < g(x) \}.\]
 
While the continuous-time process is perhaps
most natural to describe,
it is more convenient to construct the discrete-time
process first, and then to construct the continuous-time
process from that.
Thus we define a
discrete-time Markov chain,
informally obtained
by recording the locations of the successive
hits of the particle
on the boundary $\fr$.
This is a  Markov chain with state-space $\fr \cup \{ \underline
 \infty \}$ that we 
denote by
 $\xi = (\xi_n)_{n \in \Z^+}$ ($\Z^+:=\{0,1,2,\ldots\}$),
where for $\xi_n \in \fr$ we write in coordinates $\xi_n = (\xi_n^{(1)}, \xi_n^{(2)}) \in \fr$. 
Thus when $\xi_n^{(1)} > A$,
$\xi_n^{(2)}= \pm g(\xi_n^{(1)})$. We call $\xi$ the {\em collisions process}.
Informally, the state $\underline \infty$ represents the absorbing
state achieved if the particle
attains a trajectory that will never intersect with $\fr$ again.

\begin{figure}
\centering
\includegraphics[width=14cm]{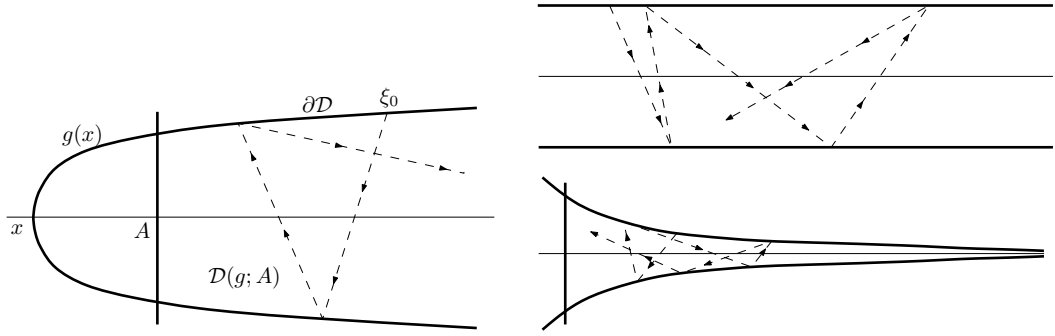}
\caption{Some examples of the tubes and the trajectories of the process}
\label{tubes}
\end{figure}

 Figure \ref{tubes} illustrates the idea of the construction.
We construct $\xi$ formally as follows.
Suppose that $\xi_0 \in \fr$, with $\xi^{(1)}_0 > 2 A$.
Let $\alpha_0, \alpha_1,\alpha_2,\ldots$
be i.i.d.~random variables
with the distribution  of a random
angle $\alpha$ where $\Pr[ |\alpha| < \pi/2 ] =1$.
 For $k \in \Z^+$,
given $\xi_k \in \fr \cup\{ \underline \infty\}$, we perform a step of
our process as follows:
\begin{itemize}
\item[(i)] If $\xi_k = \underline \infty$, set $\xi_{k+1} = \underline \infty$;
\item[(ii)] 
Otherwise, 
$\alpha_k$ specifies
a ray $\Gamma_k$ starting at $\xi_k \in \fr$ with angle $\alpha_k$ to the interior
normal to $\fr$ at $\xi_k$;
we adopt
the convention that  positive values of 
the angle correspond to the right of the normal; negative values correspond to the left.
\item[(iii)]
If the ray $\Gamma_k$ does not intersect with $\fr \setminus \{\xi_k\}$,
set $\xi_{k+1}=\underline \infty$. Otherwise,
let $(x_k,y_k)$ be the first point of intersection of the ray $\Gamma_k$ 
with $\fr$. Then if $x_k=A$, set $\xi_{k+1} = (2A,g(2A)) \in \fr$,
else set $\xi_{k+1} = (x_k,y_k) \in \fr$.
\end{itemize}
This defines the discrete-time process $\xi$.
The randomness is introduced through the random draws
from the reflection distribution $\alpha$; if
$\alpha$ is degenerate (i.e., equal to a
constant almost surely), then $\xi$ is deterministic.
Note that the construction ensures that $\xi_k^{(1)}> A$
whenever it is defined.

We need some modified form   of `reflection' away from the vertical
boundary at $x =A$,   such as that
specified by (iii). This is for technical reasons,
and in particular ensures that the process
does not jump directly to $\underline \infty$,
since that case will not be
 of interest
to us here. 
 The particular form of this
`reflection' given by (iii) is rather arbitrary; any comparable
rule will leave the behaviour unchanged.  Moreover, if the distribution
of $\alpha$ has no atom at $0$, we could instead take the reflection
rule to be the same on all boundaries without
changing the characteristics of the process.

We obtain the continuous-time stochastic process $X=(X_t)_{t \geq 0}$
with state-space $\DD \cup \fr \cup \{ \underline \infty \}$, which we
call the {\em stochastic billiard process}, from the collisions process $\xi$, essentially
by interpolation, as follows.
Suppose $X_0 = \xi_0 \in \fr$. Define
the successive {\em collision times} of the particle
by $\nu_0 :=0$, and if $\xi_k \neq \underline \infty$
\begin{equation}
\label{nudef}
 \nu_k := \sum_{j=0}^{k-1} \| \xi_{j+1} - \xi_{j} \| ,~~~ (k \in \N),
\end{equation}
where $\| \cdot\|$ is
 Euclidean distance on $\R^2$ and $\N := \{1,2,3,\ldots\}$;
if $\xi_k = \underline \infty$ set $\nu_k:= \infty$.
Then for $t \geq 0$ define
\begin{equation}
\label{ntdef}
n(t) := \max \{n \in \Z^+ : \nu_n \leq t \},
\end{equation}
so that $\nu_{n(t)} \leq t < \nu_{n(t)+1}$ for any $t \geq 0$.
Then for any $t \in (0,\infty)$ we define $X_t$ as follows:
\begin{align*}
X_t := 
\begin{cases}
\xi_n & \textrm{ if } t = \nu_n, ~n \in \N; \\
\xi_{n(t)} + \frac{ \xi_{n(t)+1} - \xi_{n(t)} }{
\| \xi_{n(t)+1} - \xi_{n(t)} \|} \cdot ( t - \nu_{n(t)}) &
\textrm{ if } t \in ( \nu_{n(t)}, \nu_{n(t)+1} ) \textrm{ and } \xi_{n(t)+1} \neq \underline \infty; \\
\underline \infty &
\textrm{ if } t \in ( \nu_{n(t)}, \nu_{n(t)+1} ) \textrm{ and } \xi_{n(t)+1} = \underline \infty.
\end{cases} 
\end{align*}
 This defines the process $X$, and in particular
$\xi$ is embedded in $X$ via $\xi_n = X_{\nu_n}$,
$n \in \Z^+$ (where $X_\infty:=\underline \infty$). 
When $X_t \in \DD \cup \fr$,
write
$X_t = (X^{(1)}_t,X^{(2)}_t)$ in coordinates. 

A realization of the sequence $(\alpha_0,\alpha_1,\ldots)$
therefore specifies the processes $\xi$ and $X$ via the construction
just described. We now describe some further assumptions that we 
make on the function $g$ that specifies the region $\DD$ and 
on the distribution of the angle $\alpha$.

Suppose that the random variable $\alpha$
for the angle of reflection
 satisfies for some $\alpha_0 \in (0,\pi/2)$
 \begin{equation}
 \label{alpha}
 \Pr [ | \alpha | < \alpha_0 ] = 1, \text{ and  }
 \Pr [ \alpha \geq x ] = \Pr [\alpha \leq -x ] \quad ( x \in \R ),
\end{equation}
 so that $\alpha$
 is bounded
 strictly away from $\pm \pi/2$ and
the distribution of $\alpha$ is symmetric around $0$.
Of special nature is the degenerate case where $\Pr [ \alpha = 0 ]=1$;
this leads to a deterministic evolution  (i.e., 
reflection
always occurs along the normals).

We will shortly introduce some assumptions on the function $g$ defining
the domain $\DD$.
If $g(x) = x$, $\DD$ is a right-cone, and
at any reflection from $\fr$
 there is positive probability (in fact,
$\Pr [ \alpha \geq 0] \geq 1/2$ under condition (\ref{alpha}))  that
the process $\xi$ (hence $X$) will go to $\underline \infty$: this is an obvious
case of transience. 
Moreover, if $\alpha$ is degenerate, then it is
clear that for $g$ strictly monotone
 $\xi_n^{(1)} \to \infty$ if and only if $g(x) \to \infty$;
so transience is possible even if $g(x)$ grows arbitrarily slowly, at
least in this deterministic case.
On the other
hand, if $g(x)=c \in (0,\infty)$, we have that if $\alpha$ is non-degenerate
then $\xi_n^{(1)}$ performs  a one-dimensional
random walk on a half-line with zero drift and uniformly bounded jumps:
so in this case $\xi$ will be null-recurrent (loosely speaking for the moment). 
On another hand, it is at least intuitively plausible
that if $g(x)$ decreases sufficiently fast, $\xi$ will
be positive-recurrent.

In the present paper we will interpolate between these three situations;
formally we assume:
\begin{itemize}
\item[(A1)] $g:[1,\infty) \to (0,\infty)$ is monotonic,
and thrice-differentiable. Moreover, there exists $\gamma \in (-\infty,1)$
for which, as $x \to \infty$,
\begin{align*}
g(x) & = x^{\gamma +o(1)} \\
g'(x) & = [ \gamma +o(1) ] \frac{g(x)}{x} \\
g''(x) & = [ \gamma (\gamma-1) + o(1) ] \frac{g(x)}{x^2}\\
|g'''(x) | & = o (x^{-2} ) .
\end{align*}
 \end{itemize}
Examples of functions satisfying (A1) include $x^\gamma$ or $x^\gamma (\log x)^\beta$, where $\gamma <1$.
Then in
particular we will
be interested in the two cases
where in addition to (A1) either:
 \begin{equation}
\label{h1} g(x) \to \infty \text{ and } x^{-1} g(x) \to 0 \text{ as } x \to \infty;
\end{equation}
or:
\begin{equation}
\label{h2} g(x) \to 0 \text{ as } x \to \infty.
\end{equation}
 
 Given (A1), a necessary condition for  (\ref{h1}) is that $\gamma \in [0,1)$; 
 a sufficient condition
 for (\ref{h2}) is $\gamma <0$. The conditions on the derivatives of $g$ imposed by (A1)
 are 
 natural smoothness 
 constraints that are necessary for our purposes.
 A particular case that we will be interested in is where
 for $\gamma < 1$,
 $g(x) = x^\gamma$ for $x \geq 1$. Then $g$ satisfies (A1)
 and moreover satisfies (\ref{h1}), (\ref{h2}) for
 $\gamma \in (0,1)$, $\gamma <0$ respectively.

When $x^{-1} g(x) \to 0$, elementary
geometry, using the fact that $\alpha$ is bounded strictly away 
from $\pm \pi/2$,
implies that $\xi$ (and hence $X$) will almost surely
 never
jump directly to $\underline \infty$, if $A$ is large enough (see Lemma \ref{cota_delta} below). 
Thus for the remainder
of the paper we can take the state-space of $\xi$ to be $\fr$
and that
of $X$ to be $\DD \cup \fr$.

In the next section we state our main results. Then in Section
\ref{open} we mention some open problems, and outline the
remainder of the paper.

\subsection{Main results}
\label{mainres}

Recall the parameter $A \in [1,\infty)$. Set
\[ \tau_A := \inf \{ t > 0 : X^{(1)}_t \leq 2 A \}; ~
\text{ and } \sigma_A :=  \inf \{ n \in \Z^+ : \xi^{(1)}_n \leq 2 A \}; \]
here and throughout the paper we use
the convention that $\inf \emptyset := + \infty$.  
We will say that the process $X$ is recurrent
if $\Pr [ \tau_A < \infty ]=1$ and transient otherwise;
if recurrent then it is positive- or null-recurrent
according to whether $\Exp [ \tau_A] < \infty$ or not. 
Similarly for the process $\xi$,
but with $\sigma_A$ instead of $\tau_A$.

By the construction of our processes,
we have that 
\begin{equation}
\label{times2}
 \tau_A < \infty \iff \sigma_A < \infty ;
\end{equation}
thus the classification of recurrence or
transience 
transfers directly between $\xi$ and $X$.

Since $|\alpha|<\alpha_0<\pi/2$, we have that the nonnegative random
variable $\tg^2 \alpha$ is bounded above by $\tg^2 \alpha_0 < \infty$. Thus
$\Exp [ \tg^2 \alpha ] < \infty$,
and $\Exp [ \tg^2 \alpha] = 0$ if and only if
$\Pr[\alpha=0]=1$. As an example,
in the case where for $\alpha_0 \in (0,\pi/2)$,
$\alpha$ is uniformly
distributed on the interval $(-\alpha_0,\alpha_0)$,
it is the case that $\Exp [ \tg^2 \alpha ] = \frac{\tg \alpha_0}{\alpha_0} - 1 >0$.

Our first result
covers the recurrence/transience
classification for our two processes $\xi$ and $X$
for the case of a growing tube. Positive-recurrence
is clearly ruled out in this case (as our processes
dominate the zero-drift process in a strip). So here
the first important issue
 is that
of transience versus null-recurrence. Define
\begin{equation}
\label{gammac}
\gamma_c   := \frac{\Exp[\tg^2\alpha]}{1+2\Exp[\tg^2\alpha]} ,
\end{equation}
so that $\gamma_c \in [0,1/2)$ and $\gamma_c =0$ if and only if $\Pr [ \alpha =0]=1$.

\begin{theo}
\label{gen_g} 
Suppose that the random variable
$\alpha$ satisfies (\ref{alpha}) and that
for $\gamma \in [0,1)$,
 $g$ satisfies (A1)  and (\ref{h1}).
Then there exists $A_0 \in (0,\infty)$ such that
for all $A>A_0$ and $\xi_0=X_0$ with $\xi^{(1)}_0 > 2A$:
\begin{itemize}
\item[(i)]
$\xi$, $X$ are transient if $\gamma > \gamma_c$; and
\item[(ii)]  $\xi$, $X$ are null-recurrent if
$\gamma < \gamma_c$.
 \end{itemize}
\end{theo}

In particular, since $\gamma_c <1/2$, Theorem
\ref{gen_g} says that for $\gamma \geq 1/2$, the processes
 $\xi$ and $X$ are always transient, regardless
of the distribution of $\alpha$.

A special case of Theorem \ref{gen_g} is given by $g(x) =(\log x)^K$ for some $K>0$;
then (A1) holds with $\gamma=0$ and (\ref{h1}) holds.
In this case it follows from Theorem \ref{gen_g} 
that $\xi$ is transient if $\Pr[\alpha=0]=1$,
otherwise it is null-recurrent. That is, the introduction
of {\em randomness} ensures that the process returns to a neighbourhood of the origin
infinitely often (with probability $1$), whereas in the deterministic
case the process is transient. The same remark
applies for other functions $g$ that grow
as $x^{o(1)}$.

Next we deal with the case of a narrowing tube, when $g$ is decreasing.
Now our processes are dominated
by the zero-drift random walk in the strip, so transience
is impossible. Thus recurrence is evident,
and the question of interest in this case is
which moments of $\tau_A, \sigma_A$ exist. The following result
is for the collisions process $\xi$.

\begin{theo}
\label{dim2_erg} 
Suppose $\alpha$ satisfies (\ref{alpha})  and that
for $\gamma \leq 0$,
 $g$ satisfies (A1)  and (\ref{h2}).
Then there exists $A_0 \in (0,\infty)$ such that
for all $A>A_0$ and $\xi_0$ with $\xi^{(1)}_0 >2A$:
\begin{itemize}
\item[(i)]
$\xi$ is positive-recurrent if
 $\gamma < -\Exp [ \tg^2 \alpha ]$; and
 \item[(ii)]  $\xi$  is null-recurrent if $\gamma > - \Exp [ \tg^2 \alpha ]$.
 \end{itemize}
 \end{theo}

It follows from
Lemma \ref{compare}(ii) below
that the continuous-time process $X$
is also positive-recurrent
under the conditions of
Theorem \ref{dim2_erg}(i). Again, Theorem \ref{dim2_erg}
shows that for a sub-polynomial
 function
such as $g(x) = 1/\log x$, a non-degenerate
$\alpha$ leads to null-recurrence.

Theorem \ref{dim2} below
deals with the important special
case
of $g(x)=x^\gamma$, $\gamma \in (0,1)$ or $\gamma <0$.
In particular, when $\alpha$ is non-degenerate,
it covers the two critical
cases omitted in Theorems \ref{gen_g}
and \ref{dim2_erg}; they are null-recurrent.

\begin{theo}
\label{dim2} Suppose $\alpha$ satisfies (\ref{alpha}) and is not degenerate, and
that $g(x)=x^\gamma$, $\gamma < 1$.
\begin{itemize}
\item[(i)] Suppose $\gamma = \gamma_c >0$. 
Then there exists $A_0 \in (0,\infty)$ such that
for all $A>A_0$ and $\xi_0$ with $\xi^{(1)}_0 >2A$,
$\xi$ is null-recurrent.
\item[(ii)] Suppose $\gamma = -\Exp[ \tg^2 \alpha] < 0$.
Then there exists $A_0 \in (0,\infty)$ such that
for all $A>A_0$ and $\xi_0$ with $\xi^{(1)}_0 >2A$,
$\xi$ is null-recurrent.
\end{itemize}
\end{theo}

We now turn to our results on almost-sure bounds for our processes, that is,
how far from the origin the particle will typically be
in both the discrete- and continuous-time processes. 
First we give a basic statement 
that says that such questions are non-trivial.

\begin{proposition}
\label{unbounded}
Suppose that $\alpha$ satisfies (\ref{alpha}). Suppose that
$g$ satisifies (A1) and
either (i) $g$ satisfies (\ref{h1});
or (ii) $g$ satisfies (\ref{h2}) and $\alpha$ is non-degenerate.
Then for $A$ sufficiently large and any $\xi_0=X_0$ with
 $\xi_0^{(1)} = X_0^{(1)} > 2A$,
a.s.,
\[ \limsup_{n \to \infty} \xi_n^{(1)} = +\infty , \text{ and }
 \limsup_{t \to \infty} X_t^{(1)} = +\infty
.\]
\end{proposition}

If $g$ satisfies (\ref{h2}), and
$\xi_0^{(1)}$ is not large enough in this last result, then 
it is clear that for some distributions
of $\alpha$ (with small bound $\alpha_0$) the particle may get `trapped' in a bounded
set. This case is not of interest for us here.

We want to quantify the result of Proposition \ref{unbounded}.
Again we can proceed
more generally, but our results are clearer if we
 take $g(x) = x^\gamma$ $(\gamma < 1$)
 from now on. Under the more general
 assumption (A1), we believe that
 our techniques still apply, with some
 modifications, and that the polynomial
 exponents in the theorems below remain valid.  
 
 Our first almost-sure bounds are for the
 collisions process in the case $\gamma \in (0,1)$:
  
\begin{theo}
\label{time1}
Suppose $\alpha$ satisfies (\ref{alpha}) and $g(x) = x^\gamma$
with
$\gamma \in (0,1)$. Suppose that $X_0^{(1)} = 
\xi_0^{(1)} > 2A$
for $A$ sufficiently large.
For any $\eps>0$,
a.s., for all but finitely many $n \in \Z^+$,
\begin{align}
\label{time1a}
\max_{0\le m\le n} \xi^{(1)}_m & \le  n^{\frac{1}{2(1-\gamma)}}(\log n)^{\frac{1}{2(1-\gamma)}+\eps}, \text{ and}\\
\label{time1b}
\max_{0\le m\le n} \xi^{(1)}_m & \geq  n^{\frac{1}{2(1-\gamma)}}(\log n)^{-\frac{1}{2(1-\gamma)}-\eps}.
\end{align}
Moreover, if $\gamma>\gamma_c$ (so that by Theorem
\ref{gen_g} $\xi$ is transient), then 
there exists $D \in (0,\infty)$ such that a.s.,
for all but finitely many $n \in \Z^+$,
\begin{equation}
\label{time1c}
\xi^{(1)}_n\ge n^{\frac{1}{2(1-\gamma)}} (\log n)^{-D}.
\end{equation}
\end{theo}
 
 Theorem \ref{time1} shows that we have polynomial behaviour, but that
 on the time-scale
 of the collisions process $\xi$, the speed $n^{-1} \xi_n^{(1)}$ can be
 very large for $\gamma$ close to $1$.
 The next result gives the corresponding
 bounds for the stochastic
 billiard process $X$. On this
 time-scale, the behaviour is very different,
 since the particle
can spend
 much more time between collisions.
 
\begin{theo}
\label{time}
Suppose $\alpha$ satisfies (\ref{alpha}) and
$g(x) =x^\gamma$ for
$\gamma \in (0,1)$. Suppose that $\xi_0^{(1)} > 2A$
for $A$ sufficiently large.
For any $\eps>0$,
a.s., for all  $t>0$ large enough,
\begin{align}
\label{time2b}
\sup_{0\le s\le t} X^{(1)}_s & \geq  t^{\frac{1}{2-\gamma}}(\log t)^{-\frac{1}{(1-\gamma)(2-\gamma)}-\eps}.
\end{align}
Moreover, if $\gamma>\gamma_c$ (so that by Theorem \ref{gen_g}
 $X$ is transient), then 
there exists $D \in (0,\infty)$ such that a.s.,
for all $t$ large enough,
\begin{equation}
\label{time2c}
t^{\frac{1}{2-\gamma}} (\log t)^{-D} \leq X^{(1)}_t\leq t^{\frac{1}{2-\gamma}} (\log t)^{D}.
\end{equation}
\end{theo}

In particular, in the (perhaps
most interesting) transient case,
(\ref{time2c})
shows that 
the asymptotic speed for the particle
 is zero,
i.e., $\lim_{t\to \infty} t^{-1} X_t^{(1)} = 0$ a.s.,
but the motion is super-diffusive. Moreover,
as $\gamma \uparrow 1$, we approach linear growth,
while as $\gamma \downarrow 0$, we approach diffusive growth,
as  expected in view of
 the remarks just above the statement of (A1).
(\ref{time2b}) also demonstrates
super-diffusive growth, even in the recurrent case.

Now we consider the case where $g(x)=x^\gamma$ for
$\gamma <0$. 
If $\gamma < - \Exp[\tg^2 \alpha]$, we use the notation
\[ \rho(\gamma) := \frac{\Exp [ \tg^2 \alpha]}{(1-2\gamma) \Exp[ \tg^2 \alpha] -\gamma} ,\]
so that $0 \leq \rho(\gamma) < 1/(2(1-\gamma))$.
The next result deals
with the collisions process $\xi$.

\begin{theo}
\label{time2}
Suppose $\alpha$ satisfies (\ref{alpha}) and
is non-degenerate, and 
$g(x) = x^\gamma$ for
$\gamma < 0$. Suppose that $\xi_0^{(1)} > 2A$
for $A$ sufficiently large.
\begin{itemize}
\item[(i)] Suppose that $\gamma > - \Exp [ \tg^2 \alpha]$. Then
for any $\eps>0$,
a.s., for all but finitely many $n \in \Z^+$,
(\ref{time1a}) and (\ref{time1b}) hold.
\item[(ii)] Suppose that $\gamma < -\Exp [ \tg^2 \alpha]$. Then
for any $\eps>0$,
a.s., for all but finitely many $n \in \Z^+$,
\begin{align}
\label{time3a}
\max_{0\le m\le n} \xi^{(1)}_m & \le  n^{\rho(\gamma)}(\log n)^{2\rho(\gamma)+\eps}, \text{ and}\\
\label{time3b}
\max_{0\le m\le n} \xi^{(1)}_m & \geq  n^{\rho(\gamma)}(\log n)^{-2\rho(\gamma)-\eps}.
\end{align}
\end{itemize}
\end{theo}

Now we state the corresponding result
for the stochastic billiard
process $X$.

\begin{theo}
\label{time3}
Suppose $\alpha$ satisfies (\ref{alpha}) and is non-degenerate,
and 
$g(x) = x^\gamma$ for
$\gamma < 0$. Suppose that $X_0^{(1)} = \xi_0^{(1)} > 2A$
for $A$ sufficiently large.
\begin{itemize}
\item[(i)] Suppose that $\gamma > - \Exp [ \tg^2 \alpha]$. Then
for any $\eps>0$,
a.s., for all $t>0$ large enough,
 (\ref{time2b}) holds.
\item[(ii)] Suppose that $\gamma < -\Exp [ \tg^2 \alpha]$. Then
for any $\eps>0$,
 a.s., for all $t>0$ large enough,
\begin{align}
 \label{time4b}
\sup_{0\le s\le t} X^{(1)}_s & \geq  t^{\frac{\rho(\gamma)}{1+\gamma \rho(\gamma)}}
(\log t)^{-\frac{4\gamma \rho(\gamma)^2 + 2 \rho(\gamma)}{1+\gamma \rho(\gamma)} - \eps}.
\end{align}
\end{itemize}
\end{theo}

Theorems \ref{time2}(ii) and \ref{time3}(ii) show that
even in the positive-recurrent case,
the behaviour is essentially
polynomial (i.e., heavy-tailed)
in nature. Moreover, the exponents
depend explicitly on the
reflection distribution $\alpha$,
unlike in the other cases.

\subsection{Open problems}
\label{open}

We briefly mention some open problems for the model
studied here. 

Of interest is the weak (distributional)
limiting behaviour of the horizontal components of
$\xi$, $X$. We have from Lemma \ref{zetalem3}
and Corollary \ref{cor1} below together
with Theorem 4.1
in \cite{lamp2}
that Lamperti's
invariance principle (\ref{invar}) below will hold
with $\eta_\bullet = (\xi^{(1)}_\bullet)^{1-\gamma}$ and
$\gamma  \in ( - \Exp [ \tg^2 \alpha], 1)$, $\alpha$ non-degenerate,
provided that Lamperti's `condition (c)' from \cite{lamp2}
is satisfied. If Lamperti's
`condition (c)' were verified, 
this would then imply that
for $\gamma \in ( -\Exp [\tg^2 \alpha], 1)$ and $\alpha$ non-degenerate,
we would have the weak invariance principle as $n \to \infty$:
\[ n^{-1/2} ( \xi^{(1)} _{\lfloor n  \bullet \rfloor } )^{1-\gamma} \Rightarrow
\Upsilon_\bullet.\]
Here $(\Upsilon_t)_{t >0}$ is a diffusion process
on $[0,\infty)$ with
Kolmogorov backwards equation
\[ u_t = \frac{a}{x} u_x + \frac{b}{2} u_{xx} ,\]
 where
\[ a = 2\gamma (1-\gamma) (1 + \Exp [ \tg^2 \alpha]),
\textrm{ and } b = 4 (1-\gamma)^2 \Exp [ \tg^2 \alpha] .\]
A more challenging problem would be to determine
whether any corresponding weak limit theory holds
for the continuous-time
process $X_t^{(1)}$ (for $\gamma > -\Exp [\tg^2 \alpha]$, say).

It may also be of interest to relax some   
of our assumptions, such as those
on the reflection distribution.
For instance, one might
consider distributions for $\alpha$
with support on all of $[-\pi/2,\pi/2]$;
there are several `natural'
distributions relevant in this case
\cite{CPSV,evans,knu}. The techniques of the present
paper require the assumption
that $\alpha$ be bounded away from $\pm \pi/2$. Indeed,
if $\alpha$ can take value $\pi/2$, say, with positive probability,
the particle will eventually jump to $\underline \infty$ in the
case of a tube of nondecreasing width. On the other hand,
if $\alpha$ has distribution on $[-\pi/2,\pi/2]$
with sufficiently light tails at the endpoints
(such that in particular $\Exp [ \tg^2 \alpha] <\infty$, say),
it may be possible to modify the techniques in the present
paper to obtain similar results;
one would need to extend various results
from processes with jumps that are bounded  to  processes
with jumps
satisfying higher-order moment assumptions, for example.

The structure of the remainder of
the paper is as follows. 
In
Section \ref{speeds}
we present and prove 
general semimartingale
criteria (of some independent interest)
on almost-sure bounds
for one-dimensional stochastic processes.
In Section \ref{lamperti}
we  make use of the results
in Section \ref{speeds}, and
give some results for the so-called
Lamperti problem,
which are of interest in their own
right as well as being
a crucial element in the proofs
of our main results.
In Section \ref{billiards}
we prove the results stated in Section \ref{mainres}
on the processes $\xi, X$.

\section{Semimartingale criteria for speeds of stochastic processes}
\label{speeds}

In this section
we present and prove
  general semimartingale criteria for obtaining
 upper and
lower almost sure bounds for discrete-time
stochastic processes on the half-line.
These results will provide
some of our main tools for the study of
 processes
with asymptotically zero mean drifts in Section \ref{lamperti} below, but for the present
section we work in some generality.

Let $({\cal F}_n)_{n \in \Z^+}$ be a filtration
on a probability space $(\Omega, {\cal F}, \Pr)$. Let $Y=(Y_n)_{n \in \Z^+}$ be a discrete-time
 $({\cal F}_n)$-adapted
stochastic process taking values in $[0,\infty)$.
Suppose that $\Pr [ Y_0 = x_0] =1$ for some $x_0 \in [0,\infty)$.

The types of process to which  our  criteria can be applied
 are quite general. For instance, due to the semimartingale nature
 of the results, we do not require that $Y$ be a Markov process. This
 fact is particularly useful if $Y$ is a process of norms $\|Z_t\|$;
 even if $(Z_t)_{t \in \Z^+}$ is Markov, $(\|Z_t\|)_{t \in \Z^+}$ will not be,
 in general. We anticipate that the general results in this section
are widely applicable. Moreover, continuous-time
processes may be treated via  embedded discrete-time
processes.
 One condition that we need
for some of our results
is that
jumps of the process $Y$ are uniformly bounded {\em above}.
Elements of the proofs in the present section
 extend ideas used
in \cite{cmp} and \cite{mw2}.

Our conditions will involve
the existence of suitable functions $f$ such that the process $f(Y)$ 
satisfies an appropriate `drift' condition. 
Our criteria will involve only first-order conditions (i.e.~expectations); no variance
results are required.
Our results
will yield almost sure bounds for $\max_{0 \leq m \leq n} Y_m$ (and
hence $Y_n$)
in terms of the function $f$ and simple functions $a, v$
that control our bounds. The functions $a, v$ will belong
to classes of eventually increasing functions defined as follows.

\begin{df} We say function $a:[1,\infty) \to [1,\infty)$ satisfies condition (C1) and
$v:[1,\infty)\to [1,\infty)$ satisfies condition (C2) if:
\begin{itemize}
\item[(C1)] $a(x) \to \infty$
 as $x \to  \infty$,
 there exists $n_a \in \N$ such that
 $x \mapsto a(x)$ is increasing for all $x \geq n_a$, and
  $\sum_{n=1}^\infty
  \frac{1}{a(n)} < \infty$;
  \item[(C2)]
$v(x) \to \infty$ as $x \to  \infty$,
 there exists $n_v \in \N$ such that
 $x \mapsto v(x)$ is increasing for all $x \geq n_v$, and
  $\sum_{n=1}^\infty
  \frac{1}{n v(n)} < \infty$.
  \end{itemize}
  \end{df}
  
Note that  the summability condition in (C2) is equivalent to the
condition that  $\sum_{n=1}^\infty  \frac{1}{v(r^n)} < \infty$
for some (hence all) $r>1$.
Throughout this section
we interpret $\log x$ as $\max \{ 1, \log x\}$. With this convention,
condition (C1) is satisfied, for example,
 by $a(x)=x^{1+\eps}$ or $x (\log x)^{1+\eps}$, where $\eps>0$, and (C2) is satisfied, for example,
  by
$v (x) = (\log x)^{1+\eps}$ or $(\log x)(\log \log x)^{1+\eps}$, $\eps>0$.

The first result in this section is a submartingale  criterion for an upper bound.

\begin{theo}
\label{thm4a}
Let $(Y_n)_{n \in \Z^+}$ be a discrete-time
 $({\cal F}_n)$-adapted
stochastic process taking values in $[0,\infty)$.
Suppose that there exists $f:[0,\infty)\to[0,\infty)$ a nondecreasing
function such that,
\begin{align}
\label{jj1}
 \Exp [ f(Y_{n+1}) - f(Y_n) \mid {\cal F}_n ] \geq 0 ~\textrm{a.s.},
 \end{align}
for all $n \in \Z^+$. Also suppose that there exists $B \in (0,\infty)$
such that for all $n \in \N$
\begin{equation}
\label{jj2} \Exp [ f(Y_n) ] \leq B n .\end{equation}
Define the nondecreasing function $f^{-1}$ for $x>0$ by
\begin{equation}
\label{finv}
 f^{-1} (x) := \sup \{ y \geq 0 : f(y) < x \}.\end{equation}
Let $a$ satisfy (C1).
Then, a.s., for all but finitely many $n \in \Z^+$,
 \begin{equation}
 \label{001x}
 \max_{0 \leq m \leq n} Y_m  \leq f^{-1}
( a(2n) ). \end{equation}
\end{theo}

A sufficient condition for (\ref{jj2}) is clearly  that
there exists $B' \in (0,\infty)$ such that
\begin{equation}
\label{jj3}
\Exp [ f(Y_{n+1}) - f(Y_n) \mid {\cal F}_n ] \leq B' ~\textrm{a.s.},\end{equation}
for all $n \in \Z^+$. We next present
a variant of Theorem \ref{thm4}, which
 relaxes the condition (\ref{jj1}) at the expense of this slightly
 stronger version of (\ref{jj2}).

\begin{theo}
Theorem \ref{thm4a} holds with conditions (\ref{jj1}) and (\ref{jj2})
replaced by the lone
condition (\ref{jj3}).
\label{thm4}
\end{theo}

For the proof of Theorem~\ref{thm4}, we need the following:
\begin{lm}
\label{lem13} Let   $(Z_n)_{n \in \Z^+}$
be an $(\F_n)$-adapted process on $[0,\infty)$ 
and $z_0 \in [0,\infty)$  
such that $\Pr [ Z_0
= z_0 ] =1$ and for some $B \in (0,\infty)$ and all $n \in \Z^+$
\begin{align}
 \label{oww3}
  \Exp [ Z_{n+1}- Z_n \mid \F_n ] \leq B ~{\rm a.s.}. 
  \end{align}
Then for any $x >0$ and any $n \in \N$ 
\begin{align}
 \label{oooo}
 \Pr \left[ \max_{0 \leq m \leq n} Z_m \geq x \right]
\leq (B n +z_0) x^{-1}.
\end{align}
\end{lm}
\proof Similarly to Doob's decomposition (see e.g.~\cite{williams},
p.~120), set $W_0:=Z_0$, and for $m \in \N$ let $W_m := Z_m +
A_{m-1}^-+A_{m-2}^-+\cdots+A_1^-$, where $A_m=  \Exp [ Z_{m+1}- Z_m
\mid \F_m ]$, $A_m^-=\max\{-A_m,0\}\ge 0$, $A_m^+=\max\{A_m,0\}\ge 0$.
Then
\[
\Exp [ W_{m+1}- W_m \mid \F_m ]=\Exp [ Z_{m+1}- Z_m+A_m^{-} \mid \F_m
]=A_m+A_m^-=A_m^+\in [0,B]
\]
so that $(W_m)$ is a nonnegative $(\F_m)$-submartingale with
$W_m\ge Z_m$ for all $m$, and $\Exp [ W_n]\le W_0+Bn=z_0+Bn$.
Hence by Doob's submartingale inequality (see
e.g.~\cite{williams}, p.~137)
\[\Pr \left[ \max_{0 \leq m \leq n} Z_m \geq x \right] \leq
\Pr \left[ \max_{0 \leq m \leq n} W_m \geq x \right] \leq x^{-1}
\Exp [ W_n] \leq (B n +z_0) x^{-1} ,\]
as required. $\qed$\\

\noindent
{\bf Proof of Theorems \ref{thm4a} and \ref{thm4}.}
First we prove Theorem \ref{thm4a}.
 Since, by (\ref{jj1}), $(f(Y_n))$
 is a nonnegative submartingale,
 Doob's
submartingale inequality 
implies that, for any $n \in \N$,
 \begin{equation}
  \label{1002dz}
 \Pr \left[ \max_{0 \leq m \leq n} f(Y_m) \geq
  a(n) \right] \leq (a(n))^{-1} \Exp[ f(Y_n)  ] \leq B  n  (a(n))^{-1},
\end{equation}
using (\ref{jj2}).
Also,
 for $n\in \N$,
 \begin{equation}
 \label{1002cz} 
 \Pr \left[
\max_{0 \leq m \leq n} f(Y_m)
\geq a(n) \right]
   = \Pr \left[ f \left( \max_{0 \leq m \leq n} Y_m \right) \geq a(n) \right] ,
   \end{equation}
  since $f$ is nondecreasing. With $f^{-1}$
  defined by (\ref{finv}),
 let $E_n$ denote the event
  \[ E_n := \left\{ \max_{0 \leq m \leq n} Y_m > f^{-1} ( a(n) )\right\} . \]
  Then since
  $z > f^{-1}(r)$  implies $f(z) > r$, we obtain
   from (\ref{1002dz}) and (\ref{1002cz}) that for all $n \in \N$
 \begin{equation}
 \label{lll1}
  \Pr  [ E_n ]
  \leq \Pr \left[ f \left( \max_{0 \leq m \leq n} Y_m \right) \geq a(n) \right] 
  \leq B n  (a(n))^{-1}  .\end{equation}
  Now
  \begin{equation}
  \label{lll2}
  \sum_{\ell=0}^\infty \frac{2^\ell}{a(2^\ell)} < \infty \iff \sum_{\ell=1}^\infty \frac{1}{a(\ell)} < \infty.\end{equation}
Hence by (\ref{lll1}) and (\ref{lll2}),
 along the subsequence $n=2^\ell$ for $\ell=0,1,2,\ldots$,
   (C1) and
  the 
   Borel-Cantelli lemma imply that, a.s., the event $E_n$
  occurs only finitely often, and in particular there exists $\ell_0 <\infty$
  such that for all $\ell \geq \ell_0$
  \[ \max_{0 \leq m \leq 2^\ell} Y_m   \leq f^{-1} ( a(2^\ell) ).\]
  Every $n \in \N$ sufficiently large has $2^{\ell_n} \leq n < 2^{\ell_n+1}$ for some $\ell_n \geq \ell_0$; then,
  a.s.,
  \[
  \max_{0 \leq m \leq n} Y_m
  \leq  \max_{0 \leq m \leq 2^{\ell_n+1}} Y_m
   \leq f^{-1} ( a(2^{\ell_n+1}) ) ,\]
   for all but finitely many $n$. Now since $2^{\ell_n+1} \leq 2n$ and
   $f^{-1}$ is nondecreasing,
(\ref{001x}) follows.

To obtain Theorem \ref{thm4}, in the previous argument we
replace (\ref{1002dz})
by an application of Lemma \ref{lem13}
with  
$Z_n = f(Y_n)$ and $x = a(n)$. $\qed$\\

We now work towards obtaining a lower bound for $\max_{0 \leq  m\leq n} Y_m$.
For
$\ell \in \N$ let $\sigma_\ell$ denote the {\em first passage time} of $\ell$ for $Y$, that is
\[
\sigma_\ell := \min \{ n \in \Z^+ : Y_n \geq \ell \}.
\]
Then $\sigma_1, \sigma_2, \ldots$ is a nondecreasing sequence of
stopping times for the process $Y$; under the condition
$\limsup_{n \to \infty} Y_n = +\infty$ a.s., $\sigma_\ell < \infty$
a.s.~for every $\ell$.

For $\ell \in \N$ and all $n \in \Z^+$, let $Y^\ell_n := Y_{n \wedge \sigma_\ell}$
(where $a \wedge b := \min \{ a,b \}$), the  process
stopped at $\sigma_\ell$; then $Y^\ell_n = Y_n$ if $n < \sigma_\ell$ and
$Y^\ell_n = Y_{\sigma_\ell} \geq \ell$ for $n \geq \sigma_\ell$.
We have the following result,  which is a `reverse Foster's criterion' (compare
 Theorem 2.1.1 in \cite{FMM}).

\begin{lm}
\label{thm5a}
Let $(Y_n)_{n \in \Z^+}$ be a discrete-time $({\cal F}_n)$-adapted
 stochastic process taking values in $[0,\infty)$, 
 such that for some $b \in \N$
    \begin{equation}
    \label{bndd} 
    \Pr [ Y_{n+1} \leq  Y_n + b ] =1,
    \end{equation}
    for all $n \in \Z^+$.
   Fix $\ell \in \N$ and $Y_0 = x_0 \in [0,\ell)$.
Suppose that there exists $f:[0,\infty) \to [0,\infty)$  a nondecreasing function
 for which, for some $\eps >0$,
  \begin{equation}
  \label{oo1}
   \Exp [ f(Y^\ell_{n+1}) - f(Y^\ell_n) \mid {\cal F}_n ] \geq \eps \1_{\{ \sigma_\ell > n \}} ~\textrm{a.s.} 
   \end{equation}
  for all $n \in \Z^+$. Then
  \[ \Exp [ \sigma_\ell ] \leq \frac{1}{\eps} f (\ell + b) .\]
  \end{lm}
  \proof Let $\ell \in \N$.
   Taking expectations in (\ref{oo1}) we have for all $m \in \Z^+$
  \[ \Exp [ f(Y^\ell_{m+1}) ] - \Exp [ f(Y^\ell_m)   ] \geq \eps \Pr [ \sigma_\ell > m ].\]
  Summing for $m$ from $0$ up to $n$ we obtain
  \[  \Exp [ f(Y^\ell_{n+1}) ] -   f(Y^\ell_0)    \geq
  \eps \sum_{m=0}^{n} \Pr [ \sigma_\ell > m ].\]
  Since $f(Y^\ell_0)= f(Y_0) \geq 0$ we obtain
  \[ \Exp[ \sigma_\ell] =\lim_{n \to \infty}
  \sum_{m=0}^{n} \Pr [ \sigma_\ell > m ]
  \leq \frac{1}{\eps} \limsup_{n \to \infty} \Exp [ f(Y^\ell_{n+1}) ]
  \leq \frac{1}{\eps} f ( \ell + b),\]
  using the fact that, by (\ref{bndd}),
  $f(Y^\ell_n) \leq f( \ell + b )$ for all $n \in \Z^+$
    since $f$  is nondecreasing.  
  $\qed$\\
 
   Now we have the following lower bound:

  \begin{theo}
   \label{thm6b}
   Let $(Y_n)_{n \in \Z^+}$ be a discrete-time $({\cal F}_n)$-adapted
 stochastic process taking values in $[0,\infty)$, such that
  condition (\ref{bndd}) holds.
Suppose that there exists $f:[0,\infty) \to [0,\infty)$ a nondecreasing function and
 $\eps >0$ for which  
 \begin{equation}
   \label{oo2}
  \Exp [ f( Y_{n+1} ) - f( Y_n )  \mid {\cal F}_n ] \geq \eps ~\textrm{a.s.}
    \end{equation}
  for all $n \in \Z^+$.
Let $v$ satisfy (C2).
   For $x\geq 0$ define the
   function $r_v$ by
 \begin{equation}
 \label{fdef3}
 r_v(x) := \inf \{ y \geq 0 : \eps^{-1} v(y) f(y+b) \geq x \}.
 \end{equation}
  Then, a.s., for
all but finitely many $n \in \Z^+$,  
\begin{equation}
 \label{003xx}
 \max_{0 \leq m \leq n} Y_m  \geq r_v(n) - b. \end{equation}
  \end{theo}
    \proof
    First note that by (C2) and the fact that $f$ is nondecreasing,
    we have that $x \mapsto r_v (x)$ is nondecreasing for all $x$ sufficiently large.
    Fix $K>1$.
    By Markov's inequality,  
\begin{equation}
\label{lll}  \Pr [ \sigma_{\ell} > v(\ell) \Exp[ \sigma_{\ell}] ]  \leq (v(\ell))^{-1} .\end{equation}
Given (\ref{lll}) and (C2),
 the 
  Borel-Cantelli lemma implies that, a.s.,
$\sigma_{ \lfloor K^\ell \rfloor} > v (\lfloor K^\ell \rfloor) \Exp [ \sigma_{\lfloor K^\ell \rfloor}]$ for only finitely
many $\ell \in \Z^+$. Moreover,
given that $f$ satisfies (\ref{oo2}), we have that (\ref{oo1}) holds
    for all $\ell \in \N$ and all $n \in \Z^+$.
Then Lemma \ref{thm5a} with
(\ref{bndd})
in this case implies that
$\Exp[ \sigma_{\ell}]  \leq \eps^{-1} f(\ell+b)$ for all $\ell$. Thus
 we have that a.s., for some $\ell_0 < \infty$ and all  $\ell \geq \ell_0$,  
$\sigma_{\lfloor K^\ell \rfloor} \leq \eps^{-1} v (\lfloor K^\ell \rfloor) f(\lfloor K^\ell \rfloor+b)$.
    Hence
with the definition of $r_v$ at (\ref{fdef3}) we have,
a.s., for all  $\ell$ sufficiently large,
\begin{align}
 \label{dddxx}  r_v ( \sigma_{\lfloor K^\ell \rfloor} ) \leq r_v (
\eps^{-1} v (\lfloor K^\ell \rfloor) f(\lfloor K^\ell \rfloor +b)) \leq \lfloor K^\ell \rfloor
\leq Y_{\sigma_{\lfloor K^\ell \rfloor}}.\end{align}
Since $\Exp [ \sigma_\ell] < \infty$ 
 for all $\ell$, we have that
$\sigma_\ell < \infty$ a.s.~for all $\ell$. Moreover,
the jumps bound (\ref{bndd})
implies that, for all $\ell \geq Y_0$, $\sigma_{\ell +b} \geq 1 + \sigma_\ell$ a.s.,
so that $\lim_{\ell \to \infty} \sigma_\ell = \infty$ a.s..
Thus, a.s.,
 every $n \in \Z^+$ 
 satisfies $\sigma_{\lfloor K^{\ell_n-1} \rfloor} \leq n < \sigma_{\lfloor K^{\ell_n} \rfloor}$
 for some $\ell_n \in \N$. Then, a.s.,
  \begin{align}
 \label{ee1xx}
  \max_{0 \leq m \leq n} Y_m \geq
 \max_{0 \leq m \leq \sigma_{\lfloor K^{\ell_n-1} \rfloor }} Y_m
 \geq Y_{\sigma_{\lfloor K^{\ell_n-1} \rfloor }} \geq \lfloor K^{\ell_n-1} \rfloor
 \geq \frac{\lfloor K^{\ell_n-1} \rfloor}{\lfloor K^{\ell_n} \rfloor} ( Y_{\sigma_{\lfloor K^{\ell_n} \rfloor}} -b),\end{align}
 for all $n$ sufficiently large,
 since $Y_{\sigma_{\lfloor K^{\ell_n} \rfloor}}  \leq \lfloor K^{\ell_n} \rfloor+b$.
 Then (\ref{ee1xx}) and (\ref{dddxx})   imply that, a.s., for all but finitely many $n \in \Z^+$,
 \[ \max_{0 \leq m \leq n} Y_m \geq
 \frac{\lfloor K^{\ell_n-1} \rfloor}{\lfloor K^{\ell_n} \rfloor} 
 \left( r_v (\sigma_{\lfloor K^{\ell_n} \rfloor}) - b \right)
  \geq
 \frac{\lfloor K^{\ell_n-1} \rfloor}{\lfloor K^{\ell_n} \rfloor}
 \left( r_v (n) - b \right) ,\]
  since $\sigma_{\lfloor K^{\ell_n} \rfloor} > n$ and $r_v(n)$ is nondecreasing for $n$ sufficiently large.
 Now taking a sequence of  values for $K$
 converging down to $1$ we obtain
  (\ref{003xx}). $\qed$

\section{Almost-sure
bounds for the Lamperti problem}
\label{lamperti}

\subsection{Introduction and results}
\label{lampres}

In this section we give upper and
 lower almost-sure bounds for  
  the so-called
  Lamperti problem of a stochastic process on the half-line with mean drift
asymptotically zero. We will need these results  
for our almost-sure bounds for the stochastic  billiard model;
our results
for the Lamperti problem
 appear to be new in the generality given here,
 and in some cases even for a nearest-neighbour
 random walk.
 
 Let $(\F_n)_{n \in \Z^+}$ be a filtration
on a
 probability space $(\Omega,\F,\Pr)$.
Let $\eta = (\eta_n)_{n \in \Z^+}$ be a discrete-time
time-homogeneous
stochastic process adapted to
$(\F_n)_{n \in \Z^+}$ and
taking values in  an unbounded
subset $\S$ of $[0,\infty)$. 
  
 We suppose that jumps of $\eta$ are uniformly bounded, that is
 there exists $B \in (0,\infty)$ such that for all $n \in \Z^+$ and all
 $x \in \S$
 \begin{equation}
 \label{jumps}
  \Pr [ | \eta_{n+1} - \eta_n | > B \mid \F_n ] = 0, ~\textrm{a.s.}.
\end{equation}
  Under the jumps condition (\ref{jumps}), 
 the jump-moment functions $\mu_1:\S \to [-B,B]$ and
$\mu_2:\S \to [0,B^2]$ given for $k \in \{1,2\}$ by
\begin{equation}
\label{mus} \mu_k (x) := \Exp [ (\eta_{n+1} - \eta_n)^k   \mid \eta_n =x], \quad
(n \in \Z^+) 
\end{equation}
are well-defined;
 in particular $\mu_2(x)$ is bounded above.
 Our basic assumption for this section
will be (A2) below.
 \begin{itemize}
 \item[(A2)]  Let $\eta$ be a discrete-time
  stochastic process on $[0,\infty)$ satisfying
  (\ref{jumps}) with $\mu_1, \mu_2$ as
  given by (\ref{mus}).
  \end{itemize}

For some of the
results in this section we also assume that
there exists $v>0$ such that, for all $x \in \S$
\begin{equation}
\label{ass2}
 \mu_2 (x) \geq v .
\end{equation}
 We will sometimes 
 make the further assumption that 
 \begin{equation}
 \label{limsup}
 \Pr \left[ \limsup_{n \to\infty} \eta_n = + \infty \right] =1 .
\end{equation}

The  model that
 we concentrate on
  is a particular case of the
 so-called
Lamperti problem
(see \cite{lamp1,lamp3,lamp2,mai})
 of a stochastic process on $[0,\infty)$ with mean drift asymptotically zero: that
 is $\mu_1(x) \to 0$ as $x \to \infty$.
 Our results do not actually assume $\mu_1(x) \to 0$, but
 it is in this case (in fact, in the case $|\mu_1(x)| = O(1/x)$)
 that they are of most interest.
 As mentioned by Lamperti \cite{lamp1},
and as is also true for
the results in this section,
only the first two
moments $\mu_1, \mu_2$ of the jump
distribution are important: there is
some form of `invariance principle' at work.
We will be mainly concerned with the case that,
from the point of view of recurrence, turns out to be
 critical; that is where $x | \mu_1 (x)|$ remains bounded
away from zero and from infinity.

 In the particular case
 of a nearest-neighbour random walk, where $\eta$ is supported on $\Z^+$,
the problem reduces
to that of a simple random walk with asymptotically zero
perturbation: 
that model was  studied
by Harris \cite{harris} and by
Hodges and Rosenblatt
\cite{hr}, and is amenable to special methods
for so-called
 birth-and-death chains. Thus many results are present
in the literature for the nearest-neighbour case: 
 in Section \ref{literature} below
we briefly mention
some of these results and their relation to the results
given in this section. For the applications
in the present paper, however, we cannot use these nearest-neighbour results.
Thus we need to prove results in the general Lamperti setting.
In fact, as we will point out below,
some of the results that we give
in the present section seem to be new even for the nearest-neighbour situation.
Thus the results that
we state in this section
are of independent interest.

As well as being of interest in their own right, stochastic processes
on the half-line with mean
drift asymptotically zero are important for the study of multidimensional processes
by the method of Lyapunov-type functions (see e.g.~\cite{FMM,mai}). For example, if $(Z_n)_{n \in \Z^+}$
is a zero-drift process with bounded jumps
on $\Z^d$, $d \geq 2$, the process $(\|Z_n\|)_{n \in \Z^+}$
 is supported on the half-line and has mean drift asymptotically zero. Clearly,
the process $\|Z_n\|$  will not be nearest-neighbour; thus the 
generality of results like Lamperti's \cite{lamp1,lamp2}
and those in the present section is
 valuable.
   
  Fix $H >0$, which will need
  to be large for some of our results. 
  We say that $\eta$ is {\em recurrent}
  or {\em transient} according
  to whether the return time $\inf \{ n \in \Z^+ : \eta_n \leq 2H\}$
  is almost surely finite or not; in the former case
  we distinguish positive- and null-recurrence
  according to whether the return time has finite or infinite
  expectation.
  
  The recurrence and transience properties of $\eta$
  were studied by Lamperti, who 
  proved the following result
   (see Theorems 3.1 and 3.2 of \cite{lamp1} with Theorem 2.1 of \cite{lamp2}; also
  Theorem 3 of \cite{mai} for  a finer result). 
  Note that all of
  the conditions that we state in this section 
  involving evaluating $\mu_1(x)$, $\mu_2(x)$
 only need apply for $x \in \S$.
  
  \begin{proposition}
  \label{lmpti}
  \cite{lamp1,lamp2}
  Suppose that (A2), (\ref{ass2}), and (\ref{limsup}) hold.
  \begin{itemize}
  \item[(i)] If for all $x>2H$, 
  $2 x | \mu_1 (x) | \leq  \mu_2 (x)$, then
  $\eta$ is null-recurrent for any $\eta_0 > 2H$.
  \item[(ii)] Suppose that there exists $\delta>0$ such that  
  for all $x>2H$ 
  \begin{equation}
 \label{mu_trans}
  2x \mu_1 (x) - \mu_2 (x) > \delta.
\end{equation}
   Then $\eta$ is transient for any $\eta_0>2H$.
  \item[(iii)] Suppose that there exists $\delta>0$ such that  
  for all $x>2H$,
  $2x \mu_1 (x) + \mu_2(x) < -  \delta$. Then $\eta$ is positive-recurrent
  for any $\eta_0 > 2H$.
  \end{itemize}
  \end{proposition}
  
For our almost-sure lower bounds,
we impose an additional `reflection'
condition that ensures that we
can avoid getting trapped in a bounded
set. Condition (A3) below
is the most appropriate
way of doing this for our applications
to the stochastic billiard model, but is clearly
stronger than is necessary for the
results in the present section.
\begin{itemize}
\item[(A3)] Given $\eta_n = x >H$,
if a jump would take $\eta_{n+1} < H$
we replace the jump with $\eta_{n+1} = 2H$
instead.
\end{itemize}

 We now state our results on almost-sure bounds.
In the general setting of the present section,
the only almost-sure bound result
that we could find in the literature
is also due to Lamperti \cite{lamp3}; see
Proposition \ref{lmpti2} in our discussion
of the literature below. Lamperti's
bound is a weaker version of
our result
Theorem \ref{prop1}(i) below,
while
Theorem \ref{prop1}(ii) gives a complementary lower bound.
 
 \begin{theo}
 \label{prop1}
 Suppose that Assumption (A2) holds.
  \begin{itemize}
   \item[(i)] 
 Suppose that there exists $C \in (0,\infty)$ such that 
 for all $x \geq 0$
 \[ 2x \mu_1(x) \leq C.\]
 Then for any $\eta_0$,
 for any $\eps>0$, a.s., for all but finitely many $n\in \Z^+$,
 \[ \max_{0 \leq m \leq n} \eta_m \leq n^{1/2} (\log n)^{(1/2)+\eps}.\]
    \item[(ii)]
 Suppose that   (A3) holds, and
 that
  there exists $\delta > 0$ such that for all $x > H$
 \[ 2x \mu_1 (x) + \mu_2 (x) \geq \delta .\]
 Then for any  $\eta_0>H$,
 for any $\eps>0$, a.s., for all but finitely many $n \in \Z^+$,
 \[ \max_{0 \leq m \leq n} \eta_m \geq n^{1/2} (\log n)^{-(1/2)-\eps}.\]
 \end{itemize}
 \end{theo}
 
 In the transient case, we prove the following lower bound
 that strengthens, in some sense,
 Theorem \ref{prop1}(ii) in this case.
 
\begin{theo}
\label{Lamp_low}
Suppose that  (A2) and (A3) 
hold.
If (\ref{mu_trans}) holds for some $\delta>0$ and all $x>H$, 
then there exists $D \in (0,\infty)$ such that,
for any $\eta_0 > H$,
 a.s., for all but finitely
many $n \in \Z^+$,
 \[
\eta_n\ge n^{1/2} (\log n)^{-D}.
\]
\end{theo}
 
 We could find no reference for a result like
 Theorem \ref{Lamp_low},
 even in the nearest-neighbour case. 
 We prove Theorem \ref{Lamp_low} in Section
 \ref{lamplowprf};
 it may be possible to extract
 bounds for the exponent $D$ in the logarithmic term
 in terms of the constants $B, \delta$ by keeping
 track of the constants in our proofs.
 Note that the transient
 critical Lamperti problem for which (A2) and (A3)  hold,
 and
 \[ \delta < 2x \mu_1 (x) - \mu_2 (x) < C \]
 for some $0<\delta<C<\infty$ and all $x$ large enough,
 satisfies the conditions
 of each of the results in Theorems \ref{prop1}
 and \ref{Lamp_low}. 
 
 The following
 discussion suggests that without
 taking into account finer behaviour (such as the
 smallest
 value of $\delta>0$ for which
 (\ref{mu_trans}) is satisfied),
 one cannot expect to improve upon Theorem \ref{Lamp_low}.
 
Let $(S^d_n)_{n \in \Z^+}$
be the symmetric simple random walk on $\Z^d$, $d \geq 2$,
so that
for $x, y \in \Z^d$,
$\Pr [ S^d_{n+1} = y \mid S^d_n = x ] = (2d)^{-1}$ if and
only if $\| y-x\|=1$. 
Then elementary calculations show that
\[ 2 \| x \| \Exp [ \| S^d_{n+1} \| - \| S^d_n \| \mid S^d_n = x ]
- \Exp [ ( \| S^d_{n+1} \| - \| S^d_n \|)^2 \mid S^d_n = x ] \to 1 - \frac{2}{d} \]
as $\| x \| \to \infty$; thus the process
$(\|S^d_n\|)_{n \in \Z^+}$ is in precisely
the critical Lamperti situation,
and for $d>2$ satisfies (\ref{mu_trans}).
 For $d>2$, 
 a  classical result of Dvoretzky and Erd\H os \cite{dver}
says that for any $\eps>0$, a.s.,
\begin{equation}
 \label{delower}
 \| S^d _n \| > n^{1/2} (\log n)^{-\frac{1}{d-2}-\eps}
\end{equation}
for all but finitely many $n \in \Z^+$, and that
this bound is
sharp in that for $\eps = 0$ the inequality (\ref{delower})
fails  infinitely often with probability $1$. In particular,
it is at least informally reasonable to
argue that by letting $d \downarrow 2$ in (\ref{delower})
one might not expect to improve in general
 upon the arbitrary
logarithmic factor in Theorem \ref{Lamp_low}.
 
 Now we state results on almost-sure bounds that we will
 need for our stochastic billiard model when $g(x)=x^\gamma$
 with $\gamma <0$. Again these results seem to be new
 in the generality given here.
  
 Theorem \ref{prop6}
  below is  further
 demonstration
 of the phenomenon of
 {\em polynomial ergodicity}
 for the critical
 Lamperti problem,
 as also
 evidenced by results of
 \cite{mp}
 in the context
 of stationary measures.
 
 \begin{theo}
 \label{prop6}
 Suppose that  (A2) and (\ref{ass2}) hold.
 \begin{itemize}
 \item[(i)]
 Suppose that there exists 
 $\kappa >1$ such that for all $x$ sufficiently large 
 \begin{equation}
 \label{c4} -2\kappa \mu_2(x) + o(1) \leq
 2x \mu_1(x)  \leq -   \kappa \mu_2(x) + o((\log x)^{-1}) .
\end{equation}
  Then for any $\eta_0$, 
  any $\eps>0$, a.s., for all but finitely many $n \in \Z^+$,
 \[ \max_{0 \leq m \leq n} \eta_m \leq n^{1/(1+\kappa)} (\log n)^{(2/(1+\kappa))+\eps}.\]
  \item[(ii)] Suppose that (A3)  holds, and
 that there exists  
  $\kappa \geq 1$ 
 such that for all $x$ sufficiently large
  \begin{equation}
 \label{c5}
    2x \mu_1 (x) + \kappa \mu_2 (x) \geq o( (\log x)^{-1}) .
   \end{equation}
  Then there exists $H \in (0,\infty)$ such that for any $\eta_0>H$,  
  for any $\eps>0$, a.s., for all but finitely many $n \in \Z^+$,
 \[ \max_{0 \leq m \leq n} \eta_m \geq n^{1/(1+\kappa)} (\log n)^{-(2/(1+\kappa))-\eps}.\]
 \end{itemize}
 \end{theo}
 
Note that the $\kappa=1$ case of Theorem \ref{prop6}(ii)
gives a
 weaker form of the lower bound
 in Theorem \ref{prop1}(ii)
 under a somewhat weaker condition. 
 
Before we prove the results stated in Section \ref{lampres},
we briefly discuss
 the existing literature
related to the present section.
  This is done in Section \ref{literature} below. Then we give the 
  proofs of Theorems \ref{prop1} and \ref{prop6}
  in Section \ref{lampprf} and Theorem \ref{Lamp_low}
  in Section \ref{lamplowprf}.

\subsection{Remarks on the literature}
\label{literature}

In the general setting of the so-called Lamperti problem,
when (A2) holds, there seem to be few almost-sure bounds
known.
The next result,   due to Lamperti himself
(see
 Theorems 2.1, 2.2, 4.2, and 5.1 in \cite{lamp3}),
includes an almost-sure upper bound  in a particular
 case where the conditions of (i) or (ii) in Proposition \ref{lmpti} hold.

\begin{proposition}
\label{lmpti2}
\cite{lamp3}
Suppose that  (A2)  holds, 
and that for
any finite interval $I \subset [0,\infty)$
\begin{equation}
\label{null}
 \lim_{n \to \infty} \frac{1}{n} \sum_{m=0}^{n-1} \Pr [ \eta_m \in I ] = 0.\end{equation}
Suppose that for $a,b \in \R$
\[ \lim_{x \to \infty} \mu_2 (x) = b >0, ~~~ \lim_{x \to \infty} (x \mu_1 (x) )= a > -(b/2) .\]
Then
for any $\eta_0$,
for any $\eps>0$, a.s., for all but finitely many $n \in \Z^+$
\begin{equation}
 \label{lampbnd}
\eta_n \leq n^{(1/2)+\eps} . \end{equation}
In addition, suppose that (\ref{null}) holds uniformly in $\eta_0$. Then
as $n \to \infty$ the following invariance principle applies:
\begin{equation}
\label{invar} 
n^{-1/2} \eta_{\lfloor n \bullet \rfloor} \Rightarrow \Upsilon_\bullet ,\end{equation}
where $(\Upsilon_t)_{t>0}$ is a diffusion
process on $[0,\infty)$ with Kolmogorov backwards equation
$u_t = (a/x) u_x + (b/2) u_{xx}$ (see Section 3 in \cite{lamp3}
for details). 
\end{proposition}

Thus Theorem \ref{prop1}(i) improves
upon the bound in (\ref{lampbnd}). Note that the condition
(\ref{null}), which does not distinguish between null-recurrence
and transience,  implies (\ref{limsup}).

A special case of the Lamperti problem on $[0,\infty)$
is the case where the process $\eta$ is supported on $\Z^+$
and only nearest-neighbour jumps are allowed. This special
case has received much more attention than the general
case described in Section \ref{lampres} above; now
we
briefly describe known results in the nearest-neighbour case. 

When they exist, these nearest-neighbour
results are sharper than the general results
that we give in Section \ref{lampres}. Thus it
may be possible
to sharpen the bounds in
the  
 results in Section \ref{lampres}.

In the nearest-neighbour case, $\eta$
  is sometimes known as
 a {\em birth-and-death}
chain (or birth-and-death random walk). Precisely, suppose that 
there exists a sequence $(p_x)_{x \in \Z^+}$ with $p_x \in (0,1)$ such that
for all $x \in \N$
and $n \in \Z^+$
\[\Pr [ \eta_{n+1} = x-1 \mid \eta_n = x ] = 
1 - \Pr [ \eta_{n+1} = x+1 \mid \eta_n = x ] = p_x,
\]
with reflection from $0$ governed by
\[ \Pr [ \eta_{n+1} = 0 \mid \eta_n = 0 ]  = 1 - \Pr [ \eta_{n+1} = 1 \mid \eta_n =0 ] = p_0. \]
(In the literature, the slightly more general model where the walk is allowed
to stay in the current position also appears. This introduces no new  essential
 features however.) 

Such processes have been extensively studied in various contexts. They are often
amenable to explicit computation; one particularly fruitful approach is via orthogonal polynomials,
dating back at least to Karlin and McGregor \cite{karlin} and subsequently employed
for instance by Voit \cite{voit90,voit93}. A recent reference
is \cite{cfr},   to which we are indebted
for bringing to our attention several of the 
papers cited in this section.

The one-step mean drift of this walk is for $x \in \N$
\[ \Exp [ \eta_{n+1} - \eta_n \mid \eta_n = x ] =  1 -2p_x .\]
We are in the Lamperti situation if we assume
$ 
\lim_{x \to \infty} p_x = 1/2$.
A case of particular interest is when $p_x \in (0,1)$ satisfies
\begin{equation}
\label{nnw}
p_x = \frac{1}{2} + \frac{\kappa}{4} x^{-\alpha} + o(x^{-\alpha}),
\end{equation}
for $\alpha >0$ and $\kappa \in \R \setminus \{ 0\}$. Then
with the notation at (\ref{mus}), $\mu_1(x) = -(\kappa/2) x^{-\alpha} + o(x^{-\alpha})$
and $\mu_2(x) =1$.

In the nearest-neighbour case, partial versions of the
recurrence result Proposition \ref{lmpti}
were given in e.g.~\cite{harris,hr}. When (\ref{nnw}) holds,
Proposition \ref{lmpti} 
implies that if $\alpha >1$, $\eta$ is null-recurrent,
while if $\alpha \in (0,1)$, $\eta$ is transient,
 positive-recurrent according to  
$\kappa <0$, $\kappa >0$. In the critical
case $\alpha =1$, $\eta$ is transient
if $\kappa <-1$, null-recurrent if $|\kappa| < 1$,
and positive-recurrent if $\kappa >1$.

Almost-sure results are known for certain
null-recurrent and transient situations. The following two propositions
collect some known results of this type;
special cases were dealt with in \cite{brez,fal81,rosen,szekely}.
The 
next result
  follows from Theorem 2.11 of \cite{voit92}; it deals
  with the transient case with $\alpha \in (0,1)$.
  (There appears to be a typo in the
  statement of Theorem 2.11 in \cite{voit92}.)
\begin{proposition}
\cite{voit92}
Suppose that (\ref{nnw}) holds for $\alpha \in (0,1)$
and $\kappa < 0$.
 Then for any $\eta_0$, a.s.
\[ \lim_{n \to \infty} \frac{\eta_n}{n^{1/(1+\alpha)}} = 
\left( \frac{|\kappa|}{2} (1+\alpha) \right)^{1/(1+\alpha)} .\]
\end{proposition}
Theorem 7.1 of Lamperti \cite{lamp3} gives a corresponding
result for more general processes (in the 
manner of Section \ref{lampres}), but with
convergence in probability only.

The following iterated logarithm result
that applies in the $\kappa<0$, $\alpha=1$ case of (\ref{nnw})
follows from Theorem 4(b) of \cite{gallardo} (compare also Theorem 1.3 of \cite{voit90}).
\begin{proposition}
\cite{gallardo}
Suppose that 
there exist $c, C$ with $0<c<C<\infty$ such that for all $x$ sufficiently large
\[ \frac{1}{2} - C x^{-1} \leq p_x \leq  \frac{1}{2} - c x^{-1}. \] 
 Then for any $\eta_0$, a.s.
\[ \limsup_{n \to \infty} \frac{\eta_n}{\sqrt{2n \log \log n}} = 1 .\]
\end{proposition}
 
Apparently missing, for example, is the complete result in the case where 
(\ref{nnw}) holds with $\alpha =1$ and $\kappa \in (0,1)$,
 for which the walk is null-recurrent;
 a one-sided result is given in
 Theorem 4(a) of \cite{gallardo}.
Also, there seem to be no (sharp) existing
 results in the critical positive-recurrent case  
 where (\ref{nnw}) holds with $\alpha =1$ and $\kappa > 1$.
 Thus Theorems \ref{prop1}, \ref{Lamp_low}, and \ref{prop6} 
add to known results even in the nearest-neighbour case.

\subsection{Proofs of Theorems \ref{prop1} and \ref{prop6}}
\label{lampprf}
 
  In order to prove Theorems \ref{prop1}
 and \ref{prop6},
 we apply the general results
 on almost-sure bounds
 for discrete-time stochastic processes
 given in Section \ref{speeds}.\\
 
 \noindent
 {\bf Proof of Theorem \ref{prop1}.}   For $x \geq 0$ and $n \in \Z^+$ we have
\begin{align}
 \Exp [ \eta_{n+1}^2  - \eta_n^2  \mid \eta_n = x] 
 &= 2x \Exp [ \eta_{n+1} -\eta_n \mid \eta_n = x]
     +\Exp [ (\eta_{n+1} -\eta_n)^2 \mid \eta_n = x] \nonumber\\
& = 2x \mu_1(x) + \mu_2(x).\label{aa1}
\end{align}
 Under the conditions of part (i) of the  theorem,
 the right-hand side of (\ref{aa1}) is uniformly bounded above. Hence 
 we can apply 
Theorem \ref{thm4}, with $Y_n = \eta_n$,
 taking $f(x)=x^2$ and $a(x)=x(\log x)^{1+\eps}$. This proves part (i). 

  Under the conditions of part (ii) of the theorem,
  we have that the right-hand side of (\ref{aa1}) is strictly positive
  for all $x > H$. Under condition
  (A3) it suffices to consider the process on $[H,\infty)$. Hence we can
  apply 
Theorem \ref{thm6b}, with $Y_n = \eta_n$,
 taking $f(x)=x^2$ and $v(x)=(\log x)^{1+\eps}$. This proves part (ii). 
$\qed$\\
  
Now we prepare for the
proof of Theorem \ref{prop6}. We need two more
lemmas, to identify suitable functions $f$ with which
to apply the criteria of Section \ref{speeds}.
 
 \begin{lm}
 \label{lem5}
 Suppose that  (A2) and (\ref{ass2}) hold.
 Suppose that there exists $\kappa>1$
 such that (\ref{c4}) holds for all $x$ sufficiently large.
   Then there exists $C \in (0,\infty)$ such that for all $x \geq 0$
   \[ \Exp [ \eta_{n+1}^{1+\kappa} (\log (1+\eta_{n+1}))^{-1}
  - \eta_{n}^{1+\kappa} (\log (1+\eta_{n}))^{-1}
   \mid \eta_n = x ] \leq C.\]
  \end{lm}
  \proof
 It follows from Taylor's theorem
applied to the function $x \mapsto x^{1+\kappa} (\log (1+x))^{-1}$
 that for $|\theta(x)|=O(1)$
as $x\to\infty$
\begin{align*}
& (x+\theta(x))^{1+\kappa} (\log (1+x+\theta(x)))^{-1}
 - x^{1+\kappa} ( \log (1+x))^{-1} \\
& = \frac{x^{1+\kappa}}{\log (1+x)}
 \left[ \frac{1+\kappa}{2x^2}
( 2x \theta(x) + \kappa \theta(x)^2 ) 
-\frac{1}{2x^2 \log (1+x)} (2x \theta(x) + (2\kappa+1) \theta(x)^2 +o(1) )
\right] .
\end{align*}
Now conditioning on $\eta_n=x$ and setting
$\theta(x)=\eta_{n+1}-\eta_n$,
taking expectations in the last
displayed equation gives
\begin{align*}
& \Exp [ \eta_{n+1}^{1+\kappa} (\log (1+\eta_{n+1}))^{-1}
  - \eta_n^{1+\kappa} (\log(1+\eta_n))^{-1} \mid \eta_n =  x]\\
& = \frac{x^{1+\kappa}}{\log (1+x)}
 \left[   
\frac{1+\kappa}{2x^2}
( 2x \mu_1(x) + \kappa \mu_2 (x) ) 
-\frac{1}{2x^2 \log (1+x)} (2x \mu_1(x) + (2\kappa+1) \mu_2(x) +o(1) )
 \right] \\
& \leq \frac{x^{\kappa-1}}{(\log (1+x))^2} 
\left[ -\frac{1}{2} \mu_2 (x) + o(1) \right],
 \end{align*}
using (\ref{c4}). Then with (\ref{ass2})
this yields the result. 
$\qed$
 
 \begin{lm}
 \label{lem6}
 Suppose that (A2) and (\ref{ass2}) hold.
  Suppose that there exists 
   $\kappa \geq 1$ such that (\ref{c5}) 
holds for all $x$ sufficiently large.
    Then there exist $H \in (0,\infty)$,
     $\eps > 0$ such that for all $x > H$ 
      \[ \Exp [ \eta_{n+1}^{1+\kappa} \log (1+\eta_{n+1})
  - \eta_n^{1+\kappa}\log(1+\eta_n) \mid \eta_n =  x] \geq \eps.\]
 \end{lm}
 \proof This time, 
it follows from Taylor's theorem that for $|\theta(x)|=O(1)$
as $x\to\infty$
\begin{align*}
& (x+\theta(x))^{1+\kappa} \log (1+x+\theta(x)) - x^{1+\kappa} \log (1+x)\\
& = x^{1+\kappa}\left[ \frac{1+\kappa}{2x^2} ( \log (1+ x) )
( 2x \theta(x) + \kappa \theta(x)^2 ) 
+\frac{1}{2x^2 } (2x \theta(x) + (2\kappa+1) \theta(x)^2 + o(1))
 \right] .
\end{align*}
Now conditioning on $\eta_n=x$ and setting
$\theta(x)=\eta_{n+1}-\eta_n$,
taking expectations in the last
displayed equation gives
\begin{align*}
& \Exp [ \eta_{n+1}^{1+\kappa} \log (1+\eta_{n+1})
  - \eta_n^{1+\kappa}\log(1+\eta_n) \mid \eta_n =  x]\\
& = x^{1+\kappa}  \left[   
\frac{1+\kappa}{2x^2} (\log (1+x))
( 2x \mu_1(x) + \kappa \mu_2 (x) ) 
+\frac{1}{2x^2 } (2x \mu_1(x) + (2\kappa+1) \mu_2(x) +o(1))
  \right] \\
& \geq x^{\kappa-1} \left[ \frac{\kappa+1}{2} \mu_2 (x) + o(1) \right]
\geq \eps > 0 , \end{align*}
for all $x$ large enough,
using (\ref{c5}), (\ref{ass2}), and the fact that $\kappa \geq 1$.
  $\qed$\\
  
  \noindent
  {\bf Proof of Theorem \ref{prop6}.} By Lemma \ref{lem5}, we can
  apply Theorem \ref{thm4} with $f(x) = x^{1+\kappa} ( \log (1+x))^{-1}$,
  $a(x)=x(\log x)^{1+\eps}$, and
  $Y_n = \eta_n$ to obtain part (i) of the theorem. On the
  other hand, by Lemma
  \ref{lem6}
  we can apply Theorem
   \ref{thm6b} with $f(x) = x^{1+\kappa} \log (1+x)$,
   $v(x)=(\log x)^{1+\eps}$, and $Y_n = \eta_n$
   to obtain part (ii) of the theorem.  
  $\qed$

\subsection{Proof of Theorem~\ref{Lamp_low}}
\label{lamplowprf}

Suppose that (\ref{mu_trans}) holds for some $\delta>0$.
Our strategy for the proof of Theorem
\ref{Lamp_low} is
to construct a scale on which the process  $\eta$
is transient with mean drift that is positive and
bounded uniformly away from $0$.
Fix $\beta>0$.
Denote $I_0 := \emptyset$ and for $r \in \N$ define the intervals
\[
I_r:=[(1+\beta)^r-B,(1+\beta)^r+B],
\]
where $B$ is as in the jumps bound (\ref{jumps}); then
for $\beta$ sufficiently large, $I_1, I_2,\ldots$ do not overlap.

We   look at the process $\eta$ at the moments at which
it enters an interval $I_r$ different from the last 
one visited.  We define the random times $\ell_1, \ell_2,\ldots$ inductively as
follows.
Set $\ell_1:= \min \{ n \in \Z^+ : \eta_n \in \cup_r I_r \}$; for $k \in \N$,
if
$\eta_{\ell_{k}}\in I_r$, let
$\ell_{k+1}:=\min\{\ell>\ell_{k}:  \eta_\ell \in I_{r-1}\cup I_{r+1}\}$.
Consider the embedded process
 $\tilde \eta = (\tilde \eta_k)_{k \in \N} = (\eta_{\ell_k})_{k \in \N}$.
 The conditions of Theorem \ref{Lamp_low}
 imply those of Theorem \ref{prop1}(ii), so in particular
 (\ref{limsup}) holds. This together with
 the jumps bound (\ref{jumps}) implies that
 \[\Pr \Big[ \bigcup_{m > n} \{ \eta_m \in I_{r-1}\cup I_{r+1}\} \mid \eta_n \in I_r \Big] =1\]
  for any $r \in \N$, so that
  the process $\tilde \eta$ is well-defined and
 the random times $\ell_k$ are almost surely finite for each $k$.
 Denote   
\begin{equation}
\label{prdef}
p_r:=\Pr[\eta_{\ell_{k+1}}\in I_{r+1} \mid \eta_{\ell_k}\in I_{r}], \quad(k \in \N).
\end{equation}
The next result says that the process $\tilde \eta$ has uniformly
positive mean drift.

\begin{lm}
\label{scale-L}
Under the conditions of Theorem \ref{Lamp_low},
there exist $\beta>0$ and $\eps_0>0$ such that 
for all $r \in \N$
\[
p_r >\frac{1}{2}+\eps_0.
\] 
\end{lm}
\proof
Let $\lambda<0$. Taylor's theorem implies that as $x \to \infty$
\[
\Exp [ \eta_{n+1}^\lambda - \eta_n^\lambda \mid \eta_n = x] 
= x^{\lambda-2} \left[  \lambda x \mu_1 (x) +  \frac{\lambda(\lambda -1)}{2}
\mu_2 (x) + O (x^{-1} ) \right] ,
\]
using (\ref{jumps}) and (\ref{mus}). Then by (\ref{mu_trans})
and (\ref{jumps}) again, we have
\[ 
\Exp [ \eta_{n+1}^\lambda - \eta_n^\lambda \mid \eta_n = x] 
\leq x^{\lambda-2} \left[  ( \lambda /2)  ( \delta  + \lambda B^2 )
+ O (x^{-1} ) \right],
\]
where $B$ is the jumps bound in (\ref{jumps}).
It follows that
for $\lambda \in ( - \delta/B^2, 0)$
\[ \Exp [ \eta_{n+1}^\lambda - \eta_n^\lambda \mid \eta_n = x] \leq 0 ,\]
for all $x$ sufficiently large. Hence for the stopping
times $\ell_k$ it is also the case that
\begin{equation}
\label{superm}
 \Exp [ \eta_{\ell_{k+1}}^\lambda - \eta_{\ell_k}^\lambda \mid \eta_{\ell_k} = x] \leq 0 ,
\end{equation}
for all $x$ sufficiently large.
Then the supermartingale property (\ref{superm})
implies that for $\beta$ large enough
and all $r \in \N$
\begin{equation}
\label{bb1} \Exp [ \eta_{\ell_{k+1}}^\lambda \mid \eta_{\ell_k} \in I_r ]
\leq [(1+\beta)^r- B]^\lambda .
\end{equation}
Also, we have that  
\begin{align}
\label{bb2}
\Exp [ \eta_{\ell_{k+1}}^\lambda \mid
\eta_{\ell_k} \in I_r ]
= &~ p_r \Exp [ \eta_{\ell_{k+1}}^\lambda \mid  \eta_{\ell_k} \in I_r,
\eta_{\ell_{k+1}} \in I_{r+1} ] \nonumber\\
& + (1-p_r) \Exp [ \eta_{\ell_{k+1}}^\lambda \mid  \eta_{\ell_k} \in I_r,
\eta_{\ell_{k+1}} \in I_{r-1} ] \nonumber\\
\geq &~ p_r [ (1+\beta)^{r+1} + B] ^\lambda
+ (1-p_r) [ (1+\beta)^{r-1} + B]^\lambda .
\end{align}
Combining (\ref{bb1}) and (\ref{bb2}) we see
that one can choose $\beta>0$, $\eps_0>0$
 such that $p_r>(1/2)+\eps_0$ for all $r$.
\qed\\

Consider, for $k \in \N$,
 \begin{equation}
 \label{zkdef}
 Z_k := \sum_{r \in \N} r \1 \{ \eta_{\ell_k} \in I_r\}.
 \end{equation}
  Then the process $Z = (Z_k)_{k \in \N}$
  is a stochastic process on $\N$ with
  nearest-neighbour
 transitions  
  that
  tracks  which interval $I_r$ the embedded
 process $\tilde \eta$ is in.
 With $p_r$ as defined by (\ref{prdef}),
we have 
\begin{equation}
\label{zjumps}
 \Pr [ Z_{k+1} = r+1 \mid Z_k = r] =
1-  \Pr [ Z_{k+1} = r-1 \mid Z_k = r] = p_r .\end{equation}
  Let $\kappa^Z_r$ denote the 
 time of the first visit of $Z$ to $r \in \N$,  i.e.
 \begin{equation}
 \label{kappaz}
  \kappa^Z_r := \min \{ k \in \N : Z_k = r \}.\end{equation}
 For $r \in \N$, $s \in \Z^+$
define
\begin{equation}
\label{gammadef}
\gamma (r,s) := \Pr [ \kappa_r^Z < \infty \mid Z_0 = r+s ].
\end{equation}
 In words, $\gamma (r,s)$ is the probability that the process $Z$
 hits $r$ in finite time, given that it starts at $r+s$. 
 
 The next result will enable us to show,
 loosely speaking,
 that the process $Z$ leaves each state 
 for good only shortly after its first visit.
 
\begin{lm}
\label{lastexit}
Under the conditions of Theorem \ref{Lamp_low},
there exists $C \in (0,\infty)$ such that
\[ \sum_{r \in \N} \gamma (r, \lfloor C \log r \rfloor) 
< \infty .\]
\end{lm}
\proof
 To estimate the required hitting probability
we introduce an auxiliary process. 
Fix $C_0 \in (0,\infty)$, which we will 
eventually take to be large. Let $r \in \N$.
Define the nonnegative  
process $(W_t)_{t \in \Z^+}$ for $t \in \Z^+$
by 
\[ W_t := \exp \left\{ - C_0^{-1} (Z_t - r)   \right\} .\] 
Then we have for $n \in \N$ and $t \in \Z^+$
\begin{align}
\label{uiop1}
\Exp [ W_{t+1} - W_t \mid Z_t = n ]
 = \exp  \{ - C_0^{-1}  ( n-r)    \} 
 \Exp [ \exp  \{ - C_0^{-1} ( Z_{t+1} - Z_t ) \} - 1
 \mid Z_t = n ] .\end{align}
We will make use of the fact that there exists $C_1  \in (0,\infty)$
such that
$\exp(-x) -1 \leq -x+ C_1 x^2$ for all $x$ with $|x| \leq 1$.
Since $Z$ has nearest-neighbour jumps,
it follows that 
\begin{align}
\label{uiop2} 
& \Exp [ \exp  \{ - C_0^{-1} ( Z_{t+1} - Z_t ) \} - 1
 \mid Z_t = n ] \nonumber\\
 \leq & - C_0^{-1} \Exp [ Z_{t+1} - Z_t \mid Z_t = n ]
+C_1 C_0^{-2} \Exp [ ( Z_{t+1} - Z_t )^2 \mid Z_t = n ]
\nonumber\\
 \leq & - C_0^{-1} (2 p_n - 1) + C_1 C_0^{-2} ,\end{align}
using (\ref{zjumps}) and the fact that $Z$ has nearest-neighbour
jumps again. It then follows
from (\ref{uiop1}), (\ref{uiop2}), and Lemma \ref{scale-L}
that we can choose $C_0$ sufficiently large,
not depending on $r$,
 so that for any $n \in \N$
\begin{equation}
\label{ww}
 \Exp [ W_{t+1} - W_t \mid Z_t = n ] \leq 0 ,\end{equation}
for all $t \in \Z^+$.

Now let $s \in \Z^+$ and take
$Z_0 = r+s$.
Recall the definition of the stopping time
 $\kappa_r^Z$ from (\ref{kappaz}).
 From (\ref{ww}) we have that $(W_{t})_{t \in \Z^+}$
is a nonnegative supermartingale,
so that $W_\infty := \lim_{t \to \infty} W_{t}$ exists
a.s., and
$\Exp [ W_0 ] \geq \Exp [ W_{\kappa_r^Z} ]$. It follows that
for any $r \in \N$, $s \in \Z^+$
\[ \Exp [ W_0 ] = \exp \{ -C_0^{-1} s \} \geq \Exp [ W_{\kappa_r^Z} \1_{\{ \kappa_r^Z < \infty\}} ]
= \gamma (r,s) ,\]
using (\ref{gammadef}).
 In particular, it follows that for some $C \in (0,\infty)$ large enough
 \[ \sum_{r \in \N} \gamma (r, \lfloor C \log r \rfloor) \leq
\sum_{r \in \N} \exp \{ - C_0^{-1} \lfloor C \log r\rfloor \} 
\leq \sum_{r \in \N} r^{-2} < \infty;\]
 completing the proof of the lemma. $\qed$\\

To complete the proof of
Theorem \ref{Lamp_low},
 we need to translate our results
for the embedded process $\tilde \eta$
back to the underlying process $\eta$. \\

\noindent
{\bf Proof of Theorem \ref{Lamp_low}.}
Consider the nearest-neighbour process
$Z=(Z_k)_{k \in \N}$  on $\N$,
as defined by (\ref{zkdef}).
Recall the definition of $\kappa^Z_r$ from (\ref{kappaz}).
 Let  $\omega^Z_r$ denote the 
 time of the   {\em last} visit of $Z$ to $r \in \N$,   i.e.
 \[   \omega^Z_r := \max \{k \in \N : Z_k = r \}.\]
 Similarly for the process $\eta$ set for $x > 0$
 \[ \kappa^\eta_x := \min \{ n \in \Z^+ : \eta_n \geq x \},
  ~~~ \omega^\eta_x := \max \{n \in \Z^+ : \eta_n \leq x \}.\]
With Lemma \ref{lastexit}, the
 Borel-Cantelli lemma implies that, a.s., for only finitely
 many $r \in \N$ does $Z$ return to $r$ after visiting $r + \lfloor C \log r \rfloor$.
So, a.s.,
 for all but finitely many $r \in \N$,
 \begin{equation}
 \label{firstlast}
  \omega^Z_r \leq \kappa^Z_{r + \lfloor C \log r \rfloor} .
\end{equation}
  Set 
\[
  r(x):= \left\lfloor \frac{\log x}{\log (1+\beta)} \right\rfloor .
\]
  Observe that by definition of the process $Z$ we have that a.s.~for
  some $C \in (0,\infty)$ and all $x$ large enough
  \[ \omega^\eta_x \leq \omega^Z_{r(x)+1} \leq \kappa^Z_{r(x)+\lfloor C \log r(x) \rfloor} ,\]
  by (\ref{firstlast}). Again by the definition of $Z$,
  it now follows that a.s.~for some $C' \in (0,\infty)$
  \begin{equation}
  \label{a1}
   \omega^\eta_x \leq \kappa^Z_{r(x)+\lfloor C \log r(x) \rfloor} \leq 
   \kappa^\eta_{\lfloor x (\log x)^{C'} \rfloor} ,
\end{equation}
for all $x$ large enough.
    
The conditions of Theorem \ref{Lamp_low}
imply those of Theorem \ref{prop1}(ii).
Hence     
the lower bound in
Theorem \ref{prop1}(ii) applies,
and so for any $\eps>0$, a.s.,
for all but finitely many $n \in \Z^+$
\[ \sup\{ x \geq 0 : \kappa^\eta_x \leq n \} \geq
\max_{0 \leq m \leq n} \eta_m \geq  
  n^{1/2} (\log n)^{-(1/2)-\eps}.\]
It follows that for any $\eps>0$, a.s.,
for all but finitely many $x \in \Z^+$
\begin{equation}
\label{a2}
 \kappa_x^\eta \leq x^2 (\log x)^{1+\eps} .
\end{equation}
So by (\ref{a1}) and (\ref{a2})
there exist $C, x_0 \in (0,\infty)$ such that, a.s., for all $x \geq x_0$
\[ \omega^\eta_x \leq x^2 (\log x)^C .\]
By the transience of $\eta$, we have that a.s.~$\eta_n \geq x_0$
for all but finitely many $n \in \Z^+$. Hence a.s.,
for all
but finitely many $n \in \Z^+$ we have
\[ n \leq \omega_{\eta_n}^\eta \leq \eta_n^2 (\log \eta_n)^C
\leq \eta_n^2 (\log n)^{C'},\]
using the jumps bound (\ref{jumps}) for the final inequality.
This proves the theorem. 
$\qed$

\section{Proofs for  stochastic billiards}
\label{billiards}

\subsection{Preliminaries}

To prove our main theorems on the stochastic billiard model,
we start by studying the properties of the process between
successive collisions; i.e., the jumps of the process $\xi$.

Suppose that at time $n \in \Z^+$ 
we have $\xi_n=(\xi_n^{(1)}, \xi_n^{(2)} )=(x, \pm g(x))$
for $x >A$, and then $\xi_n$ is reflected at angle $\alpha$ to the normal.
Denote
 $\Delta(x,\alpha):=\xi_{n+1}^{(1)}-\xi_n^{(1)}$, the jump
 of the horizontal component of $\xi$. Also set $\theta:=\arctan g'(x)$,  
so $\tg \theta= g'(x)$. 

We now proceed to obtain estimates for $\Delta(x,\alpha)$ and its moments.
The next lemma gives an upper bound on $\Delta(x,\alpha)$ that follows
from the fact that for large enough $x$ our tube will be almost flat,
while $\alpha$ is bounded strictly
away from $\pm \pi/2$.

\begin{lm}
\label{cota_delta}
Let $\alpha_0 \in (0,\pi/2)$,
and suppose that $g$ satisfies (A1), and also
(\ref{h1}) or (\ref{h2}).
Then there exist $A, C \in (0,\infty)$
such that for all $x > A$
and  
all $\alpha \in (-\alpha_0,\alpha_0)$,  
\[
|\Delta(x,\alpha)|\le C g(x).
\] 
\end{lm}
\proof
Fix $\alpha_0 \in (0,\pi/2)$.
Assuming (A1), we have that
$g'(x) \to 0$ and hence
$\theta \to 0$ as $x \to \infty$;
in particular we can choose $A$ large enough
so that for all $x   \geq A$, $|\theta| < \min \{ \alpha_0, (\pi/2)-\alpha_0\}$
and $\tg(\alpha_0+\theta)<c_0$ for some $c_0<\infty$.

\begin{figure}
\centering
\includegraphics{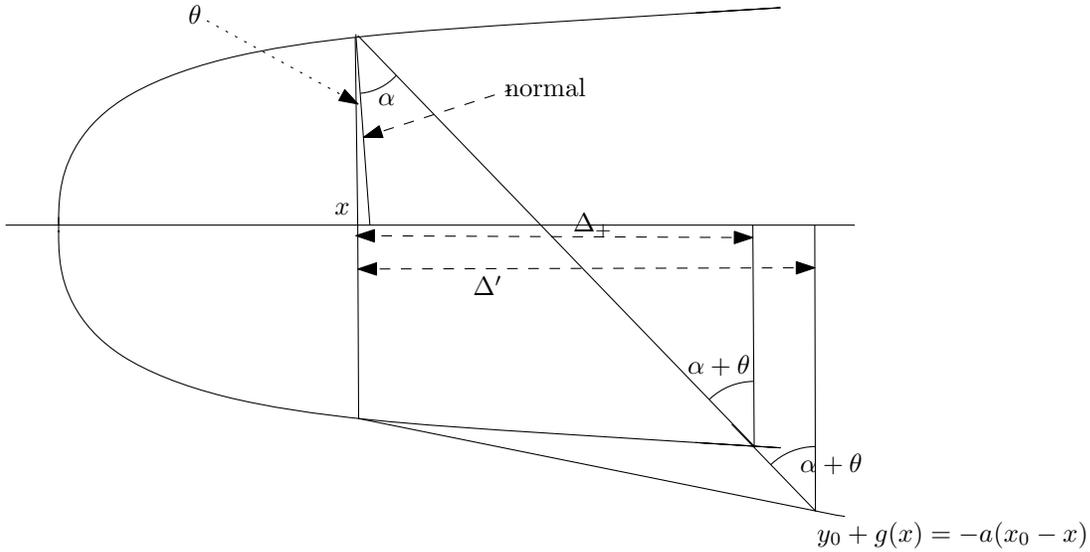}
\caption{Auxiliary construction for the proof of Lemma~\ref{cota_delta}}
\label{cota_delta_fig}
\end{figure}

By symmetry, it suffices to suppose that we start on the positive
half of the curve, i.e., at $(x,g(x))$.
We have that $|\Delta(x,\alpha)|$
is bounded by $\max\{ |\Delta(x,\alpha_0)|, |\Delta(x,-\alpha_0)|\}.$
First suppose (\ref{h1}) holds. Then $g$ is nondecreasing, so
$\theta \geq 0$.
Consider
  $\Delta(x,\alpha_0)$, represented
  as $\Delta_+$ when $\alpha=\alpha_0$
  in Figure \ref{cota_delta_fig}.
  The reflected ray at angle $\alpha_0$ to the normal
   from
  $(x,g(x))$ has equation in $(x_0,y_0)$ given by
  \[ x_0 -x = -(y_0 - g(x)) \tg (\alpha_0+\theta).\]
  
  Take $a  \in (0,1/c_0)$.
  As $g'(x)\to 0$,  we have $|g'(x)|<a$ for all $x$ large enough.
  Consider the line 
 \[ y_0 + g(x) = - a(x_0 -x).\]
This line intersects $\fr$ at  $(x,-g(x))$
and it intersects the reflected ray at angle $\alpha_0$ to the normal
from 
  $(x,g(x))$ at $x_0 \geq x$
  with
  \begin{equation}
  \label{inter1}
   x_0 - x = \frac{2g(x) \tg (\alpha_0+\theta)}
  {1-a\tg (\alpha_0+\theta)} .\end{equation}
Since $|g'(x_0)|<a$
for $x_0 \geq x$, 
the curve  $y_0 = -g(x_0)$ remains above the line  $y_0 + g(x) = - a(x_0 -x)$ for all $x_0 \geq x$.
  Thus $|\Delta(x,\alpha_0)|$
  is bounded
  by $x_0-x$ as given by (\ref{inter1});
  this is $\Delta'$ in Figure \ref{cota_delta_fig}. 
  
Thus,  for some $C \in (0,\infty)$,
$|\Delta (x,\alpha_0) | \leq  C g(x)$,
for all $x$ large enough. A similar argument
applies for $\Delta(x,-\alpha_0)$,
and with slight modification
when (\ref{h2}) holds.
\qed\\

A related geometrical argument also provides the
proof of Proposition \ref{unbounded}:\\

\noindent
{\bf Proof of Proposition \ref{unbounded}.}
Under condition (\ref{h1}),
so that $g(x)$ is strictly increasing,
the deterministic path
in the case $\Pr[\alpha=0]=1$ tends to infinity.
It is then clear that for any distribution
for $\alpha$ in this case that there exists $\eps>0$ for which
\begin{equation}
\label{updrift}
  \Pr [ \xi^{(1)}_{n+1} - \xi^{(1)}_n > \eps \mid \xi_n^{(1)} = x] > \eps 
  \end{equation}
for all $x$ large enough. The stated result for $\xi$ then follows in this case,
and hence the result for $X$ also.

Now suppose that condition (\ref{h2}) holds and 
that $\alpha$ is non-degenerate.
Then there exists $\eps_1>0$ for which $\Pr[ \alpha > \eps_1 ] > \eps_1$.
Under (A1), $g'(x) \to 0$ as $x \to \infty$, so we can
choose $A$ big enough so that the angle $\theta$ to the normal
satisfies
$|\theta| < \eps_1/2$ for all
$x > A$. From (A1), we have that $g$ is monotone and $g(x) >0$.
Then it follows that for all $x$ in any
bounded interval $(A,C)$,
(\ref{updrift})
holds for some $\eps>0$. Thus with positive
probability $\xi$ reaches any
finite horizontal distance, and the
result follows in this case also.
$\qed$ \\

The next lemma gives crucial estimates for the
first two moments of $\Delta(x,\alpha)$. These
require fairly lengthy computations.

\begin{figure}
\centering
\includegraphics{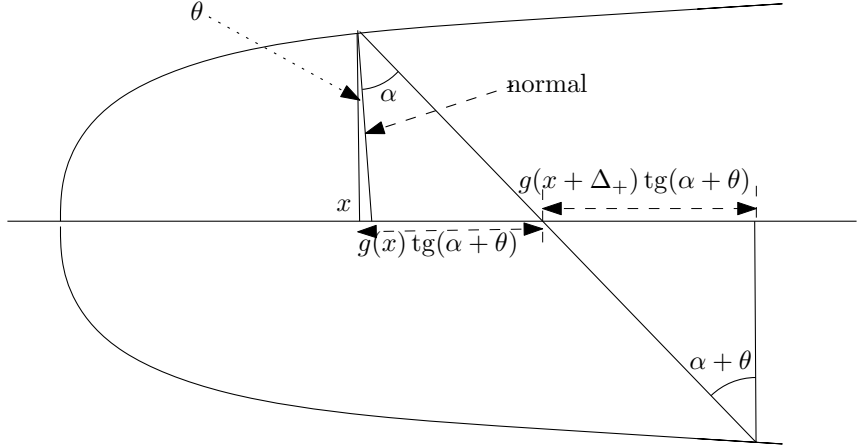}
\caption{$\Delta_+$: $\alpha > 0$}
\label{delta+}
\end{figure}

\begin{figure}
\centering
\includegraphics{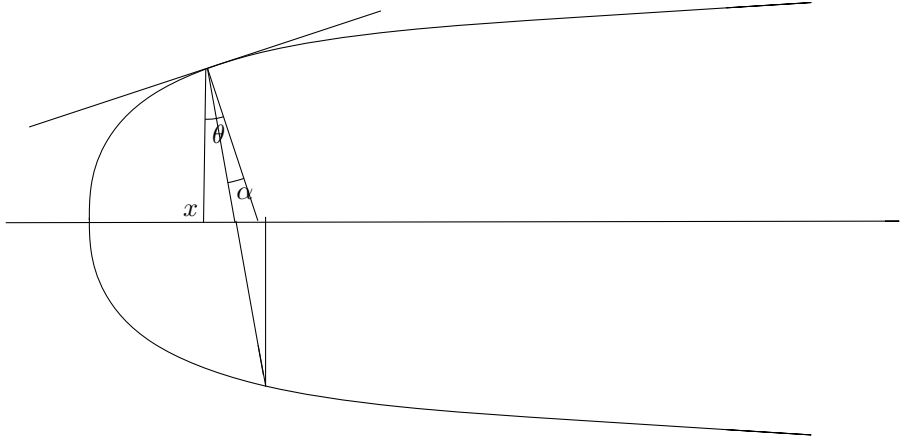}
\caption{$\Delta'_+$: $-\theta < \alpha <0$}
\label{delta'}
\end{figure}

\begin{figure}
\centering
\includegraphics{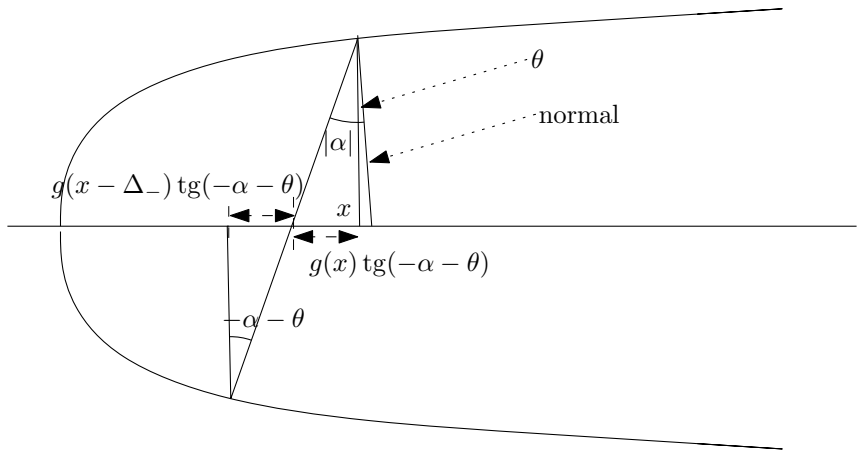}
\caption{$\Delta_-$: $\alpha < -\theta$}
\label{delta-}
\end{figure}

\begin{lm}
\label{moments_g}
 Suppose that $\alpha$ satisfies  (\ref{alpha}) and  $g(x)$ satisfies (A1),
 and also \eqref{h1} 
 or \eqref{h2}.
Then as $x \to \infty$
\begin{equation}
\label{mu1}
\Exp [\Delta(x,\alpha)]=2g'(x)g(x)(1+2\E [\tg^2\alpha]) + O(g(x)^3/x^2 ); 
\end{equation}
and
\begin{equation}
\label{mu2} 
\Exp [\Delta(x,\alpha)^2]=
      4 g(x)^2\Exp[\tg^2\alpha]+O(g(x)^3/x).
\end{equation}
\end{lm}
\proof
First suppose that (\ref{h1}) holds (the case of an increasing tube).
Then $\theta \geq 0$. Take $x$ sufficiently
large so that $\alpha_0 + \theta < \pi/2$.
Consider the jump $\Delta (x,\alpha)$.
We need to consider 3 cases.
In the first case,  $\alpha>0$ (see Figure~\ref{delta+}),
\[\Delta(x,\alpha)=\Delta_+=(g(x)+g(x+\Delta_+))\tg(\alpha+\theta).\]
In the second case, (see Figure~\ref{delta'}), $\alpha<0, |\alpha|<\theta$,
\[\Delta(x,\alpha)=\Delta'_+=(g(x)+g(x+\Delta'_+))\tg(\alpha+\theta).\]
In the third case,   (see Figure~\ref{delta-}), $\alpha<0, |\alpha|\ge\theta$,
\[-\Delta(x,\alpha)=\Delta_-=(g(x)+g(x-\Delta_-))\tg(-\alpha-\theta).\]

We start with the first of the three cases.
Using Taylor's theorem, we can write
\[
\Delta_+=\tg(\alpha+\theta)\Big[g(x)+g(x)+g'(x)\Delta_++\frac{g''(x+\phi \Delta_+)}{2}\Delta_+^2\Big],
\]
where $\phi\in[0,1]$.
Thus, writing $x_+=x+\phi \Delta_+$,
\begin{align*}
 \Delta_+&=2g(x)\frac{\tg(\alpha+\theta)}{1-g'(x)\tg(\alpha+\theta)}
            +\frac{g''(x_+)\tg(\alpha+\theta)}{2(1-g'(x)\tg(\alpha+\theta))}\Delta_+^2\\
&=2g(x)\frac{\tg\alpha+\tg \theta}{1-\tg \alpha\tg \theta
           -g'(x)(\tg\alpha+\tg \theta)}  \\
&~~+  \frac{g''(x_+)(\tg \alpha+\tg \theta)}{
         2(1-\tg \alpha\tg \theta
           -g'(x)(\tg\alpha+\tg \theta) )}\Delta_+^2\\
&=2g(x)\frac{g'(x)+\tg\alpha}{1-2g'(x)\tg \alpha
           -g'(x)^2 } + \frac{g''(x_+)}{2} 
\frac{g'(x)+\tg\alpha}{1-2g'(x)\tg \alpha
           -g'(x)^2 }\Delta_+^2,
\end{align*}
where we have used the fact that
\[
\tg(u\pm v)=\frac{\tg u\pm \tg v}{1\mp \tg u \tg v}.\]
(Recall that $\tg \theta=g'(x)$).
Analogously, in the second case,
\begin{align*}
 \Delta'_+
&=2g(x)\frac{\tg(\alpha+\theta)}{1-g'(x)\tg(\alpha+\theta)}
            +\frac{g''(x'_+)}{2}
      \frac{\tg(\alpha+\theta)}{(1-
                   g'(x)\tg(\alpha+\theta))}(\Delta'_+)^2\\
&=2g(x)\frac{g'(x)-\tg|\alpha|}{1+2g'(x)\tg |\alpha|
           -g'(x)^2 } + \frac{g''(x'_+)}{2} 
\frac{g'(x)-\tg|\alpha|}{1+2g'(x)\tg |\alpha|-g'(x)^2}(\Delta'_+)^2,
\end{align*}
where $x'_+ = x + \phi \Delta'_+$,  $\phi \in [0,1]$;
and in the third case
\begin{align*}
 \Delta_-
&=2g(x)\frac{\tg(-\alpha-\theta)}{1+g'(x)\tg(-\alpha-\theta)}
             +\frac{g''(x_-)\tg(-\alpha-\theta)}{2(1+g'(x)\tg(-\alpha-\theta))}\Delta_-^2\\
&=2g(x)\frac{-g'(x)+\tg|\alpha|}{1+2g'(x)\tg |\alpha|-g'(x)^2}+ \frac{g''(x_-)}{2} 
\frac{-g'(x)+\tg|\alpha|}{1+2g'(x)\tg |\alpha|-g'(x)^2 }\Delta_-^2,
\end{align*}
where $x_- = x - \phi \Delta_-$,  $\phi \in [0,1]$.

By Lemma~\ref{cota_delta}, we have that
$\max\{ \Delta_+ ,\Delta'_+, \Delta_- \} = O (g(x))$ a.s., and so
\[ \max \{ | x_+ - x | , | x'_+ - x | , | x_- -x | \} = O(g(x)) = o(x) .\]
So in particular $g''(x_+) = O(g''(x))$, and similarly
for $x'_+$, $x_-$.
 Thus, 
\[
\Delta_+=2g(x)\frac{g'(x)+\tg\alpha}{1-2g'(x)\tg \alpha
           -g'(x)^2} +O(g''(x)g(x)^2)
\]
\[
\Delta'_+=2g(x)\frac{g'(x)-\tg|\alpha|}{1+2g'(x)\tg |\alpha|
           -g'(x)^2}+O(g''(x)g(x)^2),
\]
and
\[
 \Delta_-= 2g(x)\frac{-g'(x)+\tg|\alpha|}{1+2g'(x)\tg |\alpha|
           -g'(x)^2} +O(g''(x)g(x)^2).
\]

Let us now estimate the first two
moments of $\Delta(x,\alpha)$. For convenience write $F(x) := \Pr [ \alpha \leq x]$.
For the first moment we obtain, using the symmetry of $F$,
\begin{align*}
\E [ \Delta(x,\alpha) ]&=\int\limits_0^{\alpha_0}\Delta_+ \ud F(\alpha) 
         +\int\limits_{-\theta}^0\Delta'_+ \ud F(\alpha) 
          -\int\limits_{-\alpha_0}^{-\theta}\Delta_- \ud F(\alpha) \\
&=\int\limits_0^{\theta} 2g(x)
      \Big[ \frac{g'(x)+\tg\alpha}{1-2g'(x)\tg \alpha
           -g'(x)^2}+ 
  \frac{g'(x)-\tg\alpha}{1+2g'(x)\tg \alpha
           -g'(x)^2}\Big]\ud F(\alpha) \\
&~~+\int\limits_{\theta}^{\alpha_0} 2 g(x)
  \Big[ \frac{g'(x)+\tg\alpha}{1-2g'(x)\tg \alpha
           -g'(x)^2}-
\frac{-g'(x)+\tg\alpha}{1+2g'(x)\tg \alpha
           -g'(x)^2}\Big]\ud F(\alpha) \\
&~~+O(g''(x)g(x)^2)\\
&=\int\limits_{0}^{\alpha_0} 2g(x)
  \Big[ \frac{g'(x)+\tg\alpha}{1-2g'(x)\tg \alpha -g'(x)^2}+
\frac{g'(x)-\tg\alpha}{1+2g'(x)\tg \alpha  -g'(x)^2}\Big]\ud F(\alpha) \\
&~~+O(g''(x)g(x)^2)\\
&=\int\limits_{0}^{\alpha_0} 2g(x)
 \Big[\frac{ 2g'(x)+4g'(x)\tg^2\alpha -2g'(x)^3}{( 1-2g'(x)\tg \alpha -g'(x)^2)
           (1+2g'(x)\tg \alpha -g'(x)^2)}  
\Big]\ud F(\alpha)\\
&~~+O(g''(x)g(x)^2). \end{align*}
The denominator in the term in square brackets in the last
integrand is $1 + O (g'(x)^2)$. Hence
\[ \E [ \Delta(x,\alpha) ]
 =2g(x)g'(x)\Big(1+2\E[\tg^2\alpha]\Big)
+O(g''(x)g(x)^2+g'(x)^3g(x)) .\]
Then using (A1) to bound the error terms, we obtain (\ref{mu1}). 

For the second moment, we have, in a similar fashion,
\begin{align*}
\E [\Delta(x, \alpha)^2]
&=\int\limits_0^{\alpha_0}\Delta_+^2 \ud F(\alpha) 
         +\int\limits_{-\theta}^0(\Delta'_+)^2 \ud F(\alpha) 
          +\int\limits_{-\alpha_0}^{-\theta}\Delta_-^2 \ud F(\alpha) \\
&=\int\limits_0^{\alpha_0}\Big[2g(x)\frac{g'(x)+\tg\alpha}{1-2g'(x)
 \tg \alpha-g'(x)^2} +O(g''(x)g(x)^2)\Big]^2\ud F(\alpha) \\
&~~+\int\limits_{-\theta}^0 \Big[2g(x)\frac{g'(x)-\tg|\alpha|}{1+2g'(x)\tg |\alpha|
           -g'(x)^2}+O(g''(x)g(x)^2)
     \Big]^2\ud F(\alpha) \\
&~~+\int\limits_{-\alpha_0}^{-\theta} \Big[2g(x)\frac{-g'(x)+\tg|\alpha|}{1+2g'(x)\tg |\alpha|
           -g'(x)^2} +O(g''(x)g(x)^2)
\Big]^2\ud F(\alpha) \\
 &=\int\limits_0^{\alpha_0}\Big[2g(x) \tg \alpha + O( g(x) g'(x) )   +O(g''(x)g(x)^2)\Big]^2\ud F(\alpha) \\
&~~+\int\limits_{-\theta}^0 \Big[2g(x) \tg \alpha +O( g(x) g'(x)) +O(g''(x)g(x)^2)
     \Big]^2\ud F(\alpha) \\
&~~+\int\limits_{-\alpha_0}^{-\theta} \Big[2g(x) (-\tg \alpha) + O(g(x) g'(x)) +O(g''(x)g(x)^2)
\Big]^2\ud F(\alpha) \\
&= 2 \int\limits_{0}^{\alpha_0} 4 g(x)^2( \tg^2\alpha ) \ud F(\alpha) 
 + O(g'(x)g(x)^2)+ O(g''(x)g(x)^3)\\
&=4 g(x)^2\E [\tg^2\alpha] + O(g'(x)g(x)^2)+ O(g''(x)g(x)^3).
\end{align*}
Then (\ref{mu2}) follows, again using (A1) to bound the error terms.

Now suppose that (\ref{h2}) holds (the case
of a decreasing tube). In this case the argument
follows similar lines to the previous one, and we only
sketch the details. Now $\theta \leq 0$.

Analogously to the case where (\ref{h1}) holds,  we define
$\Delta_+$, $\Delta_-$, $\Delta'_-$.
For  $\alpha>|\theta|$, we have 
\[\Delta(x,\alpha)=\Delta_+=(g(x)+g(x+\Delta_+))\tg(\alpha-|\theta|).\]
In the second case,  when $0\le \alpha\le |\theta|$,  
\[-\Delta(x,\alpha)=\Delta'_-=(g(x)+g(x-\Delta'_-))\tg(|\theta|-\alpha).\]
In the third case, when $\alpha<0$,
\[-\Delta(x,\alpha)=\Delta_-=(g(x)+g(x-\Delta_-))\tg(|\alpha|+|\theta|).\]

In the same way as before, we obtain
\[
\Delta_+=2g(x)\frac{g'(x)+\tg\alpha}{1-2g'(x)\tg \alpha
           -g'(x)^2} +O(g''(x)g(x)^2),
\]
\[
\Delta'_-=-2g(x)\frac{g'(x)+\tg\alpha}{1-2g'(x)\tg \alpha
           -g'(x)^2}+O(g''(x)g(x)^2),
\]
and
\[
 \Delta_-= 2g(x)\frac{-g'(x)+\tg|\alpha|}{1+2g'(x)\tg |\alpha|
           -g'(x)^2} +O(g''(x)g(x)^2).
\]
Similar computations to before then yield
the same expressions (\ref{mu1}) and (\ref{mu2}).
Lemma~\ref{moments_g} is proved.
\qed\\

The next result will allow us to
compare the recurrence times
$\sigma_A$ and $\tau_A$.

\begin{lm}
\label{compare}
Suppose that $\alpha$
satisfies (\ref{alpha}) and
$g$ satisfies (A1).
\begin{itemize}
\item[(i)] If $g$ satisfies
(\ref{h1}) then for all $A$
sufficiently large $\tau_A \geq \sigma_A$ a.s..
\item[(ii)] If $g$ satisfies (\ref{h2})
then for all $A$ sufficiently large
$\tau_A \leq \sigma_A$ a.s..
\end{itemize}
\end{lm}
\proof First we observe that 
when $\xi_n^{(1)}=x$ is large enough (such that $|g'(x)|<\tg(\frac{\pi}{2}-\alpha_0)$), we have
$\xi_{n+1}^{(2)}\xi_n^{(2)}<0$, that is,   
successive
 collisions are a.s.~on
  different sides of the tube. Thus, for $\xi_n^{(1)} = x \geq A$, a.s.,
\begin{equation}
\label{lower}
 \| \xi_{n+1} - \xi_n \| \geq | \xi_{n+1}^{(2)} - \xi_n^{(2)} |
 \geq  g (x) .\end{equation}
Now suppose that (\ref{h1}) holds.
Then
for $x$ large enough, under (\ref{h1}), (\ref{lower})
implies $\| \xi_{n+1} - \xi_n \| \geq 1$. In other
 words, the time between collisions
 for the process $X$ is no less than that
 for the process $\xi$, and part (i) follows.

Now suppose that (\ref{h2}) holds.
By the triangle inequality, we have
\[ \| \xi_{n+1} - \xi_n \| \leq
 | \xi^{(1)}_{n+1} - \xi^{(1)}_n |
 + | \xi^{(2)}_{n+1} - \xi^{(2)}_n |.\]
Thus Lemma \ref{cota_delta}
with (\ref{h2}) implies that
for some $A,C \in (0,\infty)$,
given $\xi^{(1)}_n = x \geq A$,
\[ \| \xi_{n+1} - \xi_n \| \leq
 C g(x) \leq 1, \]
 for all $x$ large enough.
 Then part (ii) follows. $\qed$

\subsection{Proofs for recurrence classification}

Our approach to studying
the horizontal component of the
collisions process $\xi$
is to consider
a rescaled version of the   process  
in such a way that we get exactly an instance of the Lamperti problem. The 
key is to find a scale on which the process has uniformly bounded jumps.

Define the function
$h: [1,\infty) \to (0,\infty)$ via
 $h(x) := x/g(x)$. Under assumption (A1) on $g$, it follows that
\begin{align*}
h'(x) & = \frac{1}{g(x)} - \frac{xg'(x)}{g(x)^2} = [ 1-\gamma + o(1) ] \frac{1}{g(x)}, \\
h''(x) & = - \frac{2 g'(x)}{g(x)^2} - \frac{x g''(x)}{g(x)^2} + \frac{2xg'(x)^2}{g(x)^3}  
  = [ \gamma (\gamma -1) + o(1) ] \frac{1}{x g(x)} ,\\
h'''(x) & = o \Big( \frac{1}{x g(x)^2} \Big) .\end{align*}

Now for $n \in \Z^+$ set 
$\zeta_n := h(\xi_n^{(1)})$.
  The process
$\zeta = (\zeta_n)_{n \in \Z^+}$ is then covered by the
  Lamperti problem (cf Section
\ref{lamperti}),
as the following result shows.
 
\begin{lm}
\label{zetalem3} 
Suppose that (A1) holds. Suppose that
$\alpha$ satisfies (\ref{alpha}).
Then there exists $B \in (0,\infty)$ such that for all $n \in \Z^+$
and all $y \geq 0$
\begin{equation}
\label{zeta0b} \Pr [ | \zeta_{n+1} - \zeta_n | \leq B \mid \zeta_n = y ] = 1.
\end{equation}
Also, for all $n \in \Z^+$, as $y \to \infty$
\begin{align}
\label{zeta1b}
m_1(y) :=
\Exp[ \zeta_{n+1}-\zeta_n \mid \zeta_n = y ]
& = \frac{2 \gamma (1-\gamma) (1+\Exp [\tg^2 \alpha ] )}{y} + o(y^{-1});
\text{ and}\\
\label{zeta2b} 
m_2 (y):=
\Exp[(\zeta_{n+1}-\zeta_n)^2\mid \zeta_n= y ]
&=4 (1-\gamma)^2 \Exp[\tg^2\alpha] +o(1).
\end{align}
\end{lm}
\proof
Given $\xi^{(1)}_n =x$,
denote $\zeta_n = y= h(x) >0$.
If the reflection is at angle $\alpha$, we have
from Taylor's theorem that as $x \to \infty$
\begin{align}
\label{def_h2} 
\zeta_{n+1}-\zeta_n & =  h(x+\Delta(x,\alpha))-h(x) \nonumber\\
& =  h'(x)\Delta(x,\alpha)+\frac{h''(x)}{2}\Delta(x,\alpha)^2 + O ( h'''(x) \Delta(x,\alpha)^3 )\nonumber\\
& =(1-\gamma +o(1) ) \frac{\Delta(x,\alpha)}{g(x)}
 +\frac{\gamma (\gamma-1) +o(1)}{2} \frac{\Delta(x,\alpha)^2}{x g(x)} + o ( g(x)/x ),
\end{align}
using Lemma~\ref{cota_delta}.
By
 Lemma~\ref{cota_delta} 
 we have that $|\Delta(x,\alpha)| = O(g(x))$, 
 and then (\ref{zeta0b}) is immediate. 
 
Taking expectations in
(\ref{def_h2}) and 
using Lemma~\ref{moments_g},
we obtain
\begin{align*}
&\E[ \zeta_{n+1}-\zeta_n \mid \zeta_n =h(x) ]\\
&=2 (1-\gamma + o(1)) (1 +2\E[\tg^2\alpha]) g'(x) 
+(2 \gamma (\gamma-1) + o(1)) \E[\tg^2\alpha] \frac{g(x)}{x}
\\
&=\frac{2\gamma(1-\gamma) (1+ \E[\tg^2\alpha] ) g(x)}{x} + o( g(x)/x).
\end{align*}
This yields (\ref{zeta1b}).
Similarly, squaring both sides of (\ref{def_h2}) and taking expectations gives (\ref{zeta2b}).
 $\qed$\\
 
 This last result,  together with
 Proposition \ref{unbounded},
 immediately implies the following:
 
   \begin{cor}
   \label{cor1}
Suppose that $g$ satisfies (A1), 
 and that $\alpha$ satisfies (\ref{alpha}). Then
 $(\zeta_n)_{n \in \Z^+}$ is a Lamperti-type problem
 as discussed in Section \ref{lamperti}, satisfying (A2) and (A3).
 Moreover, (\ref{limsup}) holds
 if $g$ satisfies (\ref{h1}),
 and also, if
 $\alpha$ is non-degenerate, if $g$ satisfies (\ref{h2}).
 Finally, (\ref{ass2}) holds provided
 that $\alpha$ is non-degenerate.
 \end{cor}
 
 In the special case where $g(x)=x^\gamma$, $\gamma <1$,
 so that $\zeta_n = (\xi_n^{(1)})^{1-\gamma}$,
  we will need
 a more precise version of Lemma \ref{zetalem3}. This is
 Lemma \ref{zetalem} below. Not only
 will this enable us to
 deal with the critical
 case in the
 recurrence classification
  (Theorem
 \ref{dim2}), it will also
 be crucial for
 our proofs of the almost-sure
 bounds carried out in Section \ref{bounds}.
 
\begin{lm}
\label{zetalem} 
Suppose that $g(x) = x^\gamma$ where
 $\gamma <1$,
 and that $\alpha$ satisfies (\ref{alpha}).
Then for all $n \in \Z^+$,
as $y \to \infty$ 
\begin{align} 
\label{zeta1}
m_1 (y) := \Exp[ \zeta_{n+1}-\zeta_n \mid \zeta_n =y ]
&=\frac{2\gamma(1-\gamma)(1+\Exp [\tg^2 \alpha ])}{y}+o(y^{-1} (\log y)^{-1}); \\
\label{zeta2} 
m_2(y) := \Exp[(\zeta_{n+1}-\zeta_n)^2\mid \zeta_n=y ]
&=4(1-\gamma)^2\Exp[\tg^2\alpha]+o((\log y)^{-1}).
\end{align}
\end{lm}
\proof
 We can apply
 Taylor's theorem to obtain,
 conditional on $\zeta_n = y=x^{1-\gamma}>0$,
 \[
 \zeta_{n+1} -\zeta_n
 = (1-\gamma) x^{-\gamma} \Delta(x,\alpha)
 - \frac{\gamma (1-\gamma)}{2} x^{-\gamma-1} \Delta (x,\alpha)^2 + O (x^{2\gamma-2}).
\]
 Then taking expectations and using (\ref{mu1}) and (\ref{mu2})
 we obtain (\ref{zeta1}). Similarly we obtain (\ref{zeta2})
 after squaring the last displayed expression. $\qed$\\
 
\noindent
{\bf Proof of Theorem \ref{gen_g}.}
Under the conditions of Theorem \ref{gen_g}, Corollary
\ref{cor1} holds.
We apply Lamperti's result
 Proposition \ref{lmpti}
to the process $\zeta$ described by Lemma \ref{zetalem3},
noting that $\zeta$ is null-recurrent, positive-recurrent, or transient exactly
when
$\xi$ is.
From Lemma \ref{zetalem3}, we have that
\begin{equation}
\label{f2}
 2 y m_1(y) + m_2 (y) = 4 (1-\gamma) ( \gamma + \Exp [\tg^2 \alpha] + o(1) ) ,\end{equation}
and also
\begin{align} 
\label{f1}
2 y m_1 (y) - m_2 (y) & = 4 (1-\gamma) \left( \gamma (1+2 \Exp [ \tg^2 \alpha]) - \Exp [ \tg^2 \alpha]
+ o(1) \right)\nonumber\\
& = 4 (1-\gamma) \left( (\gamma - \gamma_c) (1 + 2 \Exp [ \tg^2 \alpha ] ) + o(1) \right) ,\end{align}
where $\gamma_c$ is given by (\ref{gammac}).

For part (i) of the theorem, if $\gamma \in ( \gamma_c,1)$ we have from (\ref{f1})
that there exists $\delta >0$ such that
$2y m_1 (y) -m_2 (y) \geq \delta$
for all $y$ sufficiently large,
and hence
by Proposition \ref{lmpti}(ii),
$\zeta$ is transient.

For part (ii) of the theorem, it suffices to consider
the case where $\alpha$ is non-degenerate, so $\gamma_c >0$.
Then from (\ref{f1}) we have that
  for $0 \leq \gamma < \gamma_c$, $2 y m_1 (y) \leq m_2 (y)$
for all $y$ sufficiently
large. Also, in this case $\gamma + 
  \Exp [ \tg^2 \alpha] >0$, so 
we have from (\ref{f2})
that $2 y m_1 (y) \geq - m_2(y)$ for all $y$
sufficiently large.
 Hence by Proposition \ref{lmpti}(i),
$\zeta$ is null-recurrent.
 
This proves Theorem~\ref{gen_g} for the process $\xi$,
and the statement for the
process $X$ follows from 
(\ref{times2}) and Lemma \ref{compare}(i). $\qed$\\

\noindent
{\bf Proof of Theorem~\ref{dim2_erg}.} 
For part (i), if $\gamma + \Exp [ \tg^2 \alpha] <0$, we have from (\ref{f2})
 that
$2y m_1 (y) + m_2 (y) < -\delta$ for some $\delta>0$ and all $y$ sufficiently large.
Then (noting Corollary \ref{cor1})
 it follows from
Proposition \ref{lmpti}(iii)   that $\zeta$ and hence $\xi$
is positive-recurrent.

For part (ii), it suffices to suppose
that $\alpha$ is non-degenerate.
Then the $\gamma \leq 0$ case of (\ref{f1}) implies
that $2y m_1 (y) \leq m_2 (y)$ for all $y$ large enough.
On the other hand, if  
$\gamma + \Exp [ \tg^2 \alpha] >0$,
we have from (\ref{f2})
that $2 y m_1 (y) +m_2 (y) \geq 0$ for all $y$
sufficiently large. Then 
null-recurrence follows
from Proposition \ref{lmpti}(i).
\qed\\

In order to complete the proof of Theorem \ref{dim2},
we need a sharper form of Lamperti's
recurrence classification result
presented in Proposition \ref{lmpti}. Fine results
in this direction are given in \cite{mai}. We will only need the
following consequence of Theorem 3 of \cite{mai}.

\begin{lm} 
\label{mailem}
\cite{mai}
For $\eta$ a Lamperti-type problem
satisfying (A2), (\ref{ass2}), and (\ref{limsup}),
$\eta$ is null-recurrent if, for all $x$ large enough,
\[ 2x | \mu_1 (x) | \leq \left( 1 + \frac{1}{\log x} \right) \mu_2 (x) .\]
\end{lm}

\noindent{\bf Proof of Theorem \ref{dim2}.}
This now follows from Lemma \ref{mailem}
with Lemma \ref{zetalem}. $\qed$

\subsection{Proofs for almost-sure bounds}
\label{bounds}

The key to the proof of our almost-sure bounds
for the stochastic billiard model
is to
apply
our almost-sure bound results
from Section \ref{lamperti}
to  
 the rescaled
 process $\zeta$ that we studied in Lemma \ref{zetalem}.
This will allow us to obtain
Theorems \ref{time1} and \ref{time2}.
We will then derive
the results for the continuous-time process $X$, Theorems \ref{time} and \ref{time3},
from the corresponding results for $\xi$.
 Recall that 
for $n \in \Z^+$, $\zeta_n:=(\xi^{(1)}_n)^{1-\gamma}$. \\

\noindent
{\bf Proof of Theorem \ref{time1}.} The idea here is
to apply Theorem \ref{prop1} to the process $\zeta$.
By Lemma \ref{zetalem},
we have that (\ref{f2}) holds. Then since $\gamma >0$ this implies that the
conditions of Theorem \ref{prop1}(i) and (ii) are satisfied for the
process $\zeta$ (using Corollary \ref{cor1}). 
Thus for any $\eps>0$, a.s.,
for all but finitely many $n \in \Z^+$,
\[ n^{1/2} (\log n)^{-(1/2)-\eps}
\leq \max_{0 \leq m \leq n} \zeta_m \leq
n^{1/2} (\log n)^{(1/2)+\eps},\]
and then (\ref{time1a}) and  (\ref{time1b}) follow,
since $\zeta_n=(\xi^{(1)}_n)^{1-\gamma}$. 

Also by Lemma \ref{zetalem}, we have that (\ref{f1}) holds.
If $\gamma > \gamma_c$, it follows (using
Corollary \ref{cor1}) that we can
apply Theorem \ref{Lamp_low}
with $\eta = \zeta$ to obtain that for some $D \in (0,\infty)$, a.s.,
for all but finitely many $n \in \Z^+$, 
$\zeta_n \geq n^{1/2} (\log n)^{-D}$.
Then (\ref{time1c}) follows.
$\qed$\\

\noindent
{\bf Proof of Theorem \ref{time2}.} This time
we will apply Theorems \ref{prop1} and \ref{prop6}. It follows from Lemma
\ref{zetalem}
that 
\[ 2 y m_1(y) = - \kappa m_2 (y) + o((\log y)^{-1}),\]
where
\[ \kappa = \frac{-\gamma (1+ \Exp [ \tg^2 \alpha])
}{(1-\gamma) \Exp [ \tg^2 \alpha] } .\]
Hence for $\gamma < -\Exp [ \tg^2 \alpha]$, $\kappa > 1$
and (using Corollary \ref{cor1}) the conditions of parts (i) and (ii)
of Theorem \ref{prop6} are satisfied for $\zeta$. Then part (ii)
of the theorem follows.

On the other hand, for $\gamma > - \Exp [ \tg^2 \alpha]$,
we have $\kappa < 1$ so that (using Corollary \ref{cor1}) the conditions
of parts (i) and (ii) of Theorem
\ref{prop1} are satisfied. This yields part (i) of the theorem,
in the same way as in the proof of the
corresponding results in Theorem \ref{time1}.
$\qed$\\

The following lemma will enable us to derive
our `infinitely-often' lower bounds for $X^{(1)}_t$
from bounds for $\xi^{(1)}_n$.

\begin{lm}
\label{e0}
Suppose that $\gamma <1$.
Suppose that there exist $a,b >0$
with $a \gamma >-1$, such that
for any $\eps>0$, a.s.,
for all  but finitely many $n \in \Z^+$
\begin{equation}
\label{e1}
n^a (\log n)^{-b-\eps} \leq \max_{0 \leq m \leq n} \xi_m^{(1)}
\leq n^a (\log n)^{b+\eps} .\end{equation}
Then for any $\eps>0$,
a.s., for all $t$ sufficiently large,
\[ \sup_{0 \leq s \leq t} X^{(1)}_s \geq t^{\frac{a}{1+\gamma a}} (\log t)^{-\frac{2 \gamma a b +b}{1+\gamma a} -\eps} .\]
\end{lm}
\proof 
 Recall the
definition of the collision
times $\nu_k$ from (\ref{nudef}). We have from
the triangle inequality and  Lemma
\ref{cota_delta} that for some $C \in (0,\infty)$,
for all $k \in \N$,
\[ \nu_k = \sum_{j=0}^{k-1} \| \xi_{j+1} -\xi_j \| \leq
\sum_{j=0}^{k-1} \left( | \xi^{(1)}_{j+1} -\xi^{(1)}_j | + | \xi^{(2)}_{j+1} -\xi^{(2)}_j | \right) \leq
C \sum_{j=0}^k ( \xi_j^{(1)} )^ \gamma .\]
Thus by 
the upper bound in (\ref{e1}),
 for any $\eps>0$, a.s., for some $C \in (0,\infty)$ and
 all $k \in \Z^+$,
\begin{equation}
\label{d1}
 \nu_k \leq  C k^{1+\gamma a}
(\log k)^{\gamma b +\eps} .
\end{equation}
Let $\eps>0$, and for $t > 1$, set
\[ k_\eps(t) := \left \lfloor t^{\frac{1}{1+\gamma a}} 
(\log t)^{-\frac{\gamma b}{1+\gamma a}-\eps} \right \rfloor .\]
Also recall the definition of $n(t)$ from (\ref{ntdef}).
Then by (\ref{d1}) we have that for any $\eps>0$, there exist $C \in (0,\infty)$ and
$\eps'>0$ for which,
a.s.,
for all $t$ large enough
\[ \nu_{k_\eps (t)} \leq C t (\log t)^{-\eps'} \leq t;\]
hence for any $\eps>0$,
 a.s., 
for all $t$ sufficiently large 
\begin{equation}
\label{d3} 
k_\eps(t) \leq n(t), \text{ and }
 \nu_{k_\eps(t)} \leq \nu_{n(t)} \leq t < \nu_{n(t)+1}. 
\end{equation}

Now from (\ref{d3})
we have that,
a.s., for all $t$ large enough
\begin{equation}
\label{d2}
 \sup_{0 \leq s \leq t} X^{(1)}_s \geq \max_{0 \leq m \leq k_\eps (t)}
 X^{(1)}_{\nu_{m}} = \max_{0 \leq m \leq k_\eps(t)} \xi^{(1)}_{m}.
\end{equation}
Now applying the lower bound in (\ref{e1}) 
we obtain for any $\eps>0$, a.s.,
for all $t$ large enough
\[  \sup_{0 \leq s \leq t} X^{(1)}_s \geq 
(k_{\eps} (t))^a (\log k_{\eps} (t))^{-b-\eps} \geq C
 t^{\frac{a}{1+\gamma a}} (\log t)^{-\frac{\gamma a b}{1+ \gamma a} -\eps a}
(\log t )^{-b-\eps}
, \]
  using the definition of $k_\eps(t)$. Simplifying 
leads to
the desired result. $\qed$\\

We will use Lemma \ref{e0}
in the proofs of Theorems \ref{time} and \ref{time3} below.
We will apply the lemma taking either
$a = 1/(2(1-\gamma))$ with $\gamma <1$ or $a = \rho(\gamma)$ with $\gamma<0$.
Note that for $\gamma \leq 0$,
\[ \gamma \rho (\gamma) \geq \frac{\gamma}{2(1-\gamma)} \geq -\frac{1}{2},\]
so that the hypothesis $a \gamma > -1$ in Lemma
\ref{e0} is satisfied in either case.\\

\noindent
{\bf Proof of Theorem \ref{time}.}
First of all, from Theorem \ref{time1} we have that
 (\ref{time1a}) and (\ref{time1b}) hold. Thus
we can apply the $a=b=1/(2(1-\gamma))$ case of
Lemma \ref{e0}, 
which
yields (\ref{time2b}). It remains to 
 prove (\ref{time2c}).

  By the construction
of the process,
we have that a.s., for $t \geq 0$,
\begin{equation}
\label{d4}
X_t^{(1)} \in [ \min\{ \xi^{(1)}_{n(t)} , \xi^{(1)}_{n(t)+1}\} , \max\{ \xi^{(1)}_{n(t)} , \xi^{(1)}_{n(t)+1}\} ]. 
\end{equation}
 
 Suppose that $\gamma >\gamma_c$,
so that we have transience. First we prove the lower bound
in (\ref{time2c}).
 We have from (\ref{d4}) that a.s., for all $t$ large enough,
\[ X^{(1)}_t \geq \min\{ \xi^{(1)}_{n(t)} , \xi^{(1)}_{n(t)+1} \} .\]
Hence from (\ref{time1c}) we have that for some $D\in (0,\infty)$,
a.s., for all $t$ large enough
\[ X^{(1)}_t \geq (n(t))^{\frac{1}{2(1-\gamma)}} (\log n(t))^{-D}
\geq  (k_\eps(t))^{\frac{1}{2(1-\gamma)}} (\log k_\eps(t))^{-D},\]
the final inequality using (\ref{d3})  
and the fact that the function $z \mapsto z^{\frac{1}{2(1-\gamma)}} (\log z)^{-D}$
is eventually increasing in $z$.
Now (\ref{time2c}) follows by the definition of $k_\eps(t)$.

Now we prove the upper bound in (\ref{time2c}).
From (\ref{nudef}) and (\ref{lower}) we have that $\nu_k \geq \sum_{j=0}^{k-1} (\xi_j^{(1)})^\gamma$.
Then 
using (\ref{time1c}) we have that for some $D \in (0,\infty)$, a.s.,
for all but finitely many $k \in \N$,
\begin{equation}
\label{d5}
 \nu_k \geq   k^{\frac{2-\gamma}{2(1-\gamma)}} (\log k)^{-D} .\end{equation}
Let $D>0$, and for $t >1$ set
\[ k_D'(t) := \lfloor t^{\frac{2(1-\gamma)}{2-\gamma}} ( \log t)^D \rfloor .\]
Then by (\ref{d5}), we have that for $D$ large enough, a.s.,
for all $t$
sufficiently large,
\[ \nu_{k'_D (t)} \geq t, \textrm{ and } n(t) \leq k'_D (t) .\]
Now from (\ref{d4}) and (\ref{time1a})
we have that
for some $C \in (0,\infty)$, a.s., for all $t$ sufficiently large,
\[ X^{(1)}_t \leq \xi^{(1)}_{n(t)} + \xi^{(1)}_{n(t)+1}
\leq n(t)^{\frac{1}{2(1-\gamma)}} (\log n(t))^C .\]
Now using the fact that $n(t) \leq k'_D (t)$ a.s.,
and the definition of $k'_D(t)$,
the result follows.
\qed\\

\noindent
{\bf Proof of Theorem \ref{time3}.} We apply Lemma \ref{e0}
again. For part (i), we have from
part (i) of Theorem \ref{time2}
that (\ref{time1a}) and  (\ref{time1b}) hold,
so
we can apply the $a=b=1/(2(1-\gamma))$ case of
Lemma \ref{e0} to obtain (\ref{time2b}) in this case. 

For part (ii), we have from part (ii)
of Theorem \ref{time2} that (\ref{time3a}) and (\ref{time3b})
hold. So we can apply the $2a=b=2 \rho(\gamma)$ case of
Lemma \ref{e0}, which
yields (\ref{time4b}). $\qed$

\begin{center}
{\bf Acknowldegements}
\end{center}

We are grateful to the anonymous referees for their
diligence and their
detailed comments
 on an earlier version of this paper,
 which have led to several improvements.


\begin{thebibliography}{09}

\bibitem{bab} H. Babovsky, On Knudsen flows within thin tubes,
{\em J. Statist. Phys.} {\bf 44} (1986) 865--878.

\bibitem{brez} H. Br\'ezis, W. Rosenkrantz, and B. Singer,
An extension of Khintchine's estimate for large deviations
to a class of Markov chains converging
to a singular diffusion,
{\em Comm. Pure Appl. Math.} {\bf 24} (1971) 705--726.

\bibitem{cer} C. Cercignani, The Boltzmann Equation and its Applications,
Springer-Verlag, New York, 1988.

\bibitem{cmp} F. Comets, M. Menshikov, and S. Popov,
Lyapunov functions for random walks and strings in random
environment,
{\em Ann. Probab.} {\bf 26} (1998) 1433--1445.

\bibitem{CPSV} F. Comets, S. Popov, G.M. Sch\"utz, and M. Vachkovskaia (2007)
Billiards in a general domain with random reflections.
To appear in: {\it Archive for Rational Mechanics and Analysis}.
Available at {\tt arXiv.org} as {\tt math.PR/0612799}.

\bibitem{cfr} E. Cs\'aki, A. F\"oldes, and P. R\'ev\'esz (2007)
Transient NN random walk on the line. Preprint.
Available at {\tt arXiv.org} as {\tt math.PR/0707.0734}.
 
\bibitem{dver}
A. Dvoretzky and P. Erd\H os,
Some problems on random walk in space, in:
 Proceedings of the Second Berkeley Symposium on Mathematical Statistics and Probability,
  1950, pp.~353--367. University of California Press, Berkeley and Los Angeles, 1951. 

\bibitem{eml} M.D. Esposti, G.D. Magno, and M. Lenci,
An infinite step billiard,
{\em Nonlinearity} {\bf 11} (1998) 991--1013.

\bibitem{evans} S.N. Evans, Stochastic billiards on general
tables, {\em Ann. Appl. Probab.} {\bf 11} (2001) 419--437.

\bibitem{fal81} A.M. Fal$^\prime$, Certain limit theorems
for an elementary Markov random walk,
{\em Ukrainian Math. J.} {\bf 33} (1981) 433--435,
translated from {\em Ukrain. Mat. Zh.} {\bf 33} 564--566 (in Russian).

\bibitem{FMM} G. Fayolle, V.A. Malyshev, and M.V. Menshikov,
Topics in the Constructive Theory of Countable Markov Chains,
Cambridge University Press, 1995.
 
\bibitem{gallardo} L. Gallardo,
Comportement asymptotique des marches al\'eatoires associ\'ees
aux polynomes de Gegenbauer et applications,
{\em Adv. Appl. Probab.} {\bf 16} (1984) 293--323 (in French).

\bibitem{harris} T.E. Harris,
First passage and recurrence distributions,
{\em Trans. Amer. Math. Soc.} {\bf 73} (1952) 471--486.

\bibitem{hr} J.L. Hodges, Jr. and M. Rosenblatt,
Recurrence-time moments in random walks,
{\em Pacific J. Math.} {\bf 3} (1953) 127--136.
 
\bibitem{karlin} S. Karlin and J. McGregor,
Random walks,
{\em Illinois J. Math.} {\bf 3} (1959) 66--81.

\bibitem{knu} M. Knudsen, Kinetic Theory of Gases: Some Modern
Aspects, Methuen's Monographs on Physical Subjects, Methuen, London, 1952.

\bibitem{lamp1} J. Lamperti, 
Criteria for the recurrence and
transience of stochastic processes I, 
{\em J. Math. Anal. Appl.} {\bf 1} (1960) 314--330.

\bibitem{lamp3} J. Lamperti, A new class
of probability limit theorems,
{\em J. Math. Mech.} {\bf 11} (1962) 749--772.

\bibitem{lamp2} J. Lamperti, Criteria for
stochastic processes II: passage-time moments, 
{\em J. Math. Anal. Appl.} {\bf 7} (1963) 127--145.

\bibitem{lenci1} M. Lenci,
Escape orbits for non-compact flat billiards,
{\em Chaos} {\bf 6} (1996) 428--431.

\bibitem{lenci2} M. Lenci,
Semi-dispersing billiards with an infinite 
cusp I, {\em Comm. Math. Phys.} {\bf 230}
(2002) 133--180.

\bibitem{lenci3} M. Lenci,
Semidispersing billiards with an infinite
cusp. II, {\em Chaos} {\bf 13} (2003)
105--111.

\bibitem{mai} M.V. Menshikov, I.M. Asymont, and R. Iasnogorodskii, 
Markov processes with asymptotically zero drifts, 
{\em Problems of
Information Transmission} {\bf 31} (1995) 248--261, translated from
{\em Problemy Peredachi Informatsii} {\bf 31} 60--75 (in
Russian).

\bibitem{mp} M.V. Menshikov and S.Yu. Popov, 
Exact power estimates for countable Markov chains, 
{\em Markov Processes Relat. Fields} {\bf 1} (1995) 57--78.

\bibitem{mw2} M.V. Menshikov and A.R. Wade,
Logarithmic speeds for one-dimensional perturbed
random walk in random environment,
{\em Stochastic Processes Appl.}
{\bf 118} 
 (2008) 389--416.

\bibitem{rosen} W.A. Rosenkrantz,
A local limit theorem for a certain class
of random walks,
{\em Ann. Math. Statist.} {\bf 37} (1966) 855--859.

\bibitem{szekely} G.J. Sz\'ekely,
On the asymptotic properties of diffusion processes,
{\em Ann. Univ. Sci. Budapest. E\"otv\"os Sect. Math.} {\bf 17}
(1974) 69--71.

\bibitem{tab} S. Tabachnikov,
Billiards,
Soci\'et\'e Math\'ematique de France,
Paris, 1995.

\bibitem{voit90} M. Voit,
A law of the iterated logarithm for a class of polynomial
hypergroups,
{\em Monatsh. Math.} {\bf 109} (1990) 311--326.

\bibitem{voit92} M. Voit,
Strong laws of large numbers for random walks associated with a class
of one-dimensional convolution structures,
{\em Monatsh. Math.} {\bf 113} (1992) 59--74.

\bibitem{voit93} M. Voit,
A law of the iterated logarithm for Markov chains on $\mathbb{N}_0$ associated with
orthogonal polynomials,
{\em J. Theoret. Probab.} {\bf 6} (1993) 653--669.
 
 \bibitem{williams} D. Williams,
 Probability With Martingales,
  Cambridge University Press, 1991.

\end{thebibliography}
\end{document}